\documentclass[graybox,12pt]{svmult}

\evensidemargin\oddsidemargin

\usepackage{mathptmx}       
\usepackage{helvet}         
\usepackage{courier}        

\usepackage{makeidx}         
\usepackage{graphicx}        
\usepackage{multicol}        

\usepackage{amsmath}
\usepackage{amsfonts}
\usepackage{color}     
\usepackage{url}
\usepackage{verbatim}
\usepackage{enumerate}

\makeindex             

\usepackage{amssymb}
\usepackage{amsmath}
\usepackage{epsfig}
\usepackage{mathptmx}       
\usepackage{helvet}         
\usepackage{courier}        
\usepackage{type1cm}        

\usepackage{makeidx}         
\usepackage{graphicx}        
\usepackage{multicol}        
\usepackage[bottom]{footmisc}

\usepackage{amsbsy}
\usepackage{MnSymbol}
\usepackage{amsfonts}
\usepackage{amsopn}
\usepackage{dsfont}
\usepackage[english]{babel}
\usepackage[applemac]{inputenc}
\usepackage[colorlinks=true,breaklinks=true]{hyperref}
\usepackage[T1]{fontenc}
\usepackage{xcolor}
\usepackage{enumerate}
\usepackage{latexsym,amsmath,eufrak}
\usepackage{xcolor}
\usepackage{stackrel}

\begin{document}
\title*{New monotonicity and infinite divisibility properties for the Mittag-Leffler function and for the stable distributions}
\author{Nuha Altaymani, Wissem Jedidi
\thanks{Department of Statistics \& OR, King Saud University, P.O. Box 2455, Riyadh 11451, Saudi Arabia, \email{wjedidi@ksu.edu.sa}, \email{nuha.altaymani@gmail.com}}
}

\maketitle
\date{\today}
\abstract{{\bf \!:}  Hyperbolic complete monotonicity  property ($\mathrm{HCM}$) is a way to check if a distribution is a generalized gamma  ($\mathrm{GGC}$), hence is infinitely divisible.
In this work, we illustrate to which extent the Mittag-Leffler functions $E_\alpha, \;\alpha \in (0,2]$,
enjoy the  $\mathrm{HCM}$  property, and then intervene deeply in the probabilistic context. We prove that, for suitable  $\alpha$ and complex numbers $z$,  the real and imaginary part of the functions $x\mapsto E_\alpha \big(z x\big)$, are tightly linked to the stable distributions and to the generalized Cauchy kernel. }
\smallskip

\makeatletter
\newsavebox\myboxA
\newsavebox\myboxB
\newlength\mylenA

\newcommand*\xoverline[2][2]{%
    \sbox{\myboxA}{$\m@th#2$}%
    \setbox\myboxB\null
    \ht\myboxB=\ht\myboxA%
    \dp\myboxB=\dp\myboxA%
    \wd\myboxB=#1\wd\myboxA
    \sbox\myboxB{$\m@th\overline{\copy\myboxB}$}
    \setlength\mylenA{\the\wd\myboxA}
    \addtolength\mylenA{-\the\wd\myboxB}%
    \ifdim\wd\myboxB<\wd\myboxA%
       \rlap{\hskip 0.5\mylenA\usebox\myboxB}{\usebox\myboxA}%
    \else
        \hskip -0.5\mylenA\rlap{\usebox\myboxA}{\hskip 0.5\mylenA\usebox\myboxB}%
    \fi}
\makeatother
\def\i{{\infty}}
\def \l{{\lambda}}

\def\i{{\infty}}
\def\g{{\mathbb{G}}}
\def\exs{{\mathbb{G}}_1}
\def\G{{\Gamma}}
\newcommand{\aar}{\mathrm{a}}
\newcommand{\bbr}{\mathrm{b}}
\newcommand{\ccr}{\mathrm{c}}
\newcommand{\qq}{\mathrm{q}}
\newcommand{\dd}{\mathrm{d}}
\newcommand{\cop}{\cos(\pi\al)}
\newcommand{\Nset}{\mathds{N}}
\newcommand{\Zset}{\mathds{Z}}
\newcommand{\Rset}{\mathds{R}}
\newcommand{\Cset}{\mathds{C}}
\newcommand{\oi}{(0,\infty)}
\newcommand{\RP}{\mathbb{R}_+}
\newcommand{\Ima}{\mathfrak{Im}}
\newcommand{\Rel}{\mathfrak{Re}}
\newcommand{\Lap}{\mathfrak{L}}
\newcommand{\pab}[1]{\big(#1\big)}
\newcommand{\pa}[1]{\left(#1\right)}
\newcommand{\et}{\text{\:}\text{ \normalfont and }\text{\:}}
\newcommand{\Exp}[1]{\exp\left(#1\right)}
\newcommand{\tapu}{  {\mathbb{T}_a}^{\scriptscriptstyle 1/(1+a)}  }
\newcommand{\tapm}{ {\mathbb{T}_a}^{\scriptscriptstyle 1/(1-a)} }
\newcommand{\tal}{\mathbb{T}_{\alpha}}
\newcommand{\tud}{\mathbb{T}_{1/2}}
\newcommand{\tumal}{\mathbb{T}_{1-\alpha}}
\newcommand{\talm}{\mathbb{T}_{\alpha-1}}
\newcommand{\Tt}{\mathbb{T}_{t}}
\newcommand{\taa}{\mathbb{T}_{a}}
\newcommand{\T}{\mathbb T}
\newcommand{\al}{\alpha}
\newcommand{\la}{\lambda}
\newcommand{\sta}{\mathbb{S}}
\newcommand{\stabb}{\sta_{\alpha,\beta}}
\newcommand{\sal}{\mathbb{S}_{\alpha}}
\newcommand{\sat}{\mathbb{S}_{t}}
\newcommand{\tat}{\mathbb{T}_{t}}
\newcommand{\dac}{\mathbb{D}_{c,t}}
\newcommand{\Be}{\mathbb{B}}
\newcommand{\ga}{\mathbb{G}}
\newcommand{\Ga}{\Gamma}
\newcommand{\er}{\mathds{E}}
\newcommand{\pr}{\mathds{P}}
\newcommand{\ii}{{\rm i}}
\newcommand{\jj}{{\rm j}}
\newcommand{\jc}{\jj_c}
\newcommand{\jcb}{\xoverline{\jj}_c}
\newcommand{\II}{{\rm  1~\hspace{-1.2ex}l}}
\newcommand{\simdis}{\stackrel{{\rm d}}{=}}
\newcommand{\limdis}{\stackrel{{\rm d}}{\longrightarrow}}
\newcommand{\siminfini}{\underset{+\infty}{\sim}}
\newcommand{\CM}{\mathcal{CM}}
\newcommand{\ST}{\mathcal{S}}
\newcommand{\BF}{\mathcal{BF}}
\newcommand{\CB}{\mathcal{CBF}}
\newcommand{\TB}{\mathcal{TBF}}
\newcommand{\ID}{\mathrm{I\,D}}
\newcommand{\SD}{\mathrm{SD}}
\newcommand{\BO}{\mathrm{BO}}
\newcommand{\HCM}{\mathrm{HCM}}
\newcommand{\hcm}{\mathcal{HCM}}
\newcommand{\GGC}{\mathrm{GGC}}
\newtheorem{teo}{Theorem}[section]
\newtheorem{prop}[teo]{Proposition}
\newtheorem{propr}[teo]{Property}
\newtheorem{lem}[teo]{Lemma}
\newtheorem{ex}[teo]{Example}
\newtheorem{rmk}[teo]{Remark}
\newtheorem{cor}[teo]{Corollary}
\newtheorem{defi}[teo]{Definition}
\newenvironment{proofp}{\noindent {\bf Proof. }}{\ \ \ $\square$\\}
\newenvironment{proofof}{\noindent {\bf Proof  of }}{\ \ \ $\square$\\}
\long\def\symbolfootnote[#1]#2{\begingroup
\def\thefootnote{\fnsymbol{footnote}}\footnote[#1]{#2}\endgroup}
\section{Introduction and first results}

The Mittag-Leffler function
\begin{equation}\label{mittagd}
E_a (z)= \sum_{k \geq 0} \frac{z^k}{\Gamma (1+k a)}, \quad z\in \mathbb{C}, \,a >0\,,
\end{equation}
is a widely studied special function, see \cite{main} and the references therein. The essence of this work is to unveil several probabilistic interpretation for the functions  $E_\alpha, \;\alpha \in (0,2]$. To this end we need some setting.

\subsection{Bernstein functions and infinite divisibility}
A function $f:\oi \to \oi$ is {\it completely monotone} (we denote $f\in \CM$) if it is infinitely differentiable and it satisfies $(-1)^n f^{(n)}(x)\geq 0, \;\; \mbox{for all} \; n\in \Nset, \; x>0.$
By Bernstein characterisation, $f \in \CM$ if, and only if, it is the  Laplace transform of some  measure $$f(\lambda):= \Lap_\tau(\lambda)=\int_{[0,\infty)} e^{-\lambda x} \,\tau(dx), \quad \lambda >0.$$
A subset of $\CM$ is  the class $\ST_a$ of {\it generalized Stieltjes transforms} of order $a>0$, viz.  of functions $f$ represented by
\begin{equation}\label{repsti}
f(\lambda)=\dd + \frac{\qq}{\lambda^a}+\int_{(0,\infty)} \frac{1}{(\lambda+x)^a} U(dx),\quad \lambda >0,
\end{equation}
where $\qq , \dd\geq 0$ and $U$ is a Radon measure on $(0,\infty)$ such that $\int_{\oi}(1+u)^{-a}\;U (du)<\infty.$   Stieltjes transforms of order 1, are usually called {\it Stieltjes functions}. The class $\BF$ of {\it Bernstein functions} consists of functions of the form
\begin{equation}\label{BF}
     \phi(\lambda)  =\qq+ \dd \lambda + \int_{\oi}(1-e^{-\lambda x})\,\Pi(dx)
      =\qq+ \dd \lambda + \lambda \int_{\oi} e^{-\lambda x}\,\Pi(x, \infty)\;dx,\quad \lambda\geq 0,
    \end{equation}
where $\qq\geq 0$ is the so-called {\it killing rate}, $\dd \in \mathds{R}$ is the {\it drift} and $\Pi$ is the {\it L\'evy measure} of $\phi$, i.e. a positive measure on $\oi$ which satisfies $\int_{\oi}\left(x \wedge 1\right) \Pi(dx) <\infty$. Two important subclasses of $\BF$ are the ones of {\it Thorin Bernstein functions} and of {\it complete Bernstein functions}, denoted by $\TB$ and $\CB$, respectively. A function $\phi$ belongs to $\TB$ (respectively $\CB$)  if it is represented  by
\begin{equation}\label{repphi3tt}
\phi(\lambda)= \qq+ \dd \,\lambda + \,\int_{\oi}  \log \left(1+\frac{\lambda }{u}\right) \;U(du),\quad
\big(\mbox{respectively}\;\; \phi(\lambda)= \qq+ \dd  \, \lambda + \,\int_{\oi}  \frac{\lambda }{\lambda +u} \; V (du)\big), \quad \lambda\geq 0,
\end{equation}
where the positive measures $U$ and $V$ satisfy
\begin{equation}\label{lolo}
\int_{(0,1)}|\log u|\,U(du)+ \int_{[1,\infty)} u^{-1}\;U(du) <\infty, \; \;\mbox{and} \;\;\frac{V(du)}{u} \;\mbox{is  a L\'evy measure.}
\end{equation}
Note the equivalences
\begin{equation}\label{eqstielcbt}
\phi \in \BF \; \mbox{(respectively $\TB, \;\CB $)} \Longleftrightarrow \phi \geq 0 \;\mbox{and $\phi' \in \CM$  (respectively $\ST_1, \;\ST_2$)}.
\end{equation}
The book of Hirschman \& Widder \cite{HW} is a good reference for the classes $\ST_a$ and, especially for $\ST_1$, we refer to \cite[Chapter 2]{SSV}. The books of Schilling,  Song \& Vondra\c cek \cite{SSV} and the one of Steutel \& van Harn \cite{steutel} are our main references for the class $\BF$ and its subclasses, and for  {\it infinitely divisible measures} as well. The distribution of a nonnegative random variable $X$ is said  to be  infinitely divisible, and we denote $X\sim \ID$, if there exists an  i.i.d. sequence $(X_i^n)_{1\leq i\leq n}$, such that $X\simdis X_1^n+\ldots +X_n^n$. The celebrated {\it L\'evy-Khintchine formula} gives the following characterization through the cumulant function of $X\geq 0$:
\begin{equation}\label{bb}
X\sim \ID \;  \Longleftrightarrow \;\er[e^{-\lambda X}]^t=\left(\int_{[0,\infty)}e^{- \lambda x} \, \pr(X\in dx)\right)^t \in \CM, \;\; \forall t>0 \; \Longleftrightarrow \; \phi_X(\lambda) :=-\log \er[e^{-\lambda X}]  \in \BF.
\end{equation}
\subsection{The $\HCM$ property and $\GGC$ distributions}
The class of infinitely divisible distributions behind $\TB$ (respectively $\CB$) is known as the {\it generalized gamma convolutions} $\GGC$ (respectively the {\it Bondesson class} $\BO$), and  it corresponds to the smallest class of sub-probability measures on $[0,\infty)$ which contains mixtures of gamma  (respectively exponential) distributions and which is closed under convolutions and vague limits. See \cite[Theorem 9.7 and  Proposition 9.11]{SSV}. These classes were introduced by Olaf Thorin,  and were widely developed by Lennart Bondesson in \cite{B}, see also \cite{SSV,steutel}. We will now introduce  an important subclass of $\GGC$. A function $f: (0,+\infty)\to(0,+\infty)$ is said to be {\it hyperbolically completely monotone}, and we denote $f\in \hcm$ if, for every $u > 0$, the function
$$H_u(w):= f(uv)f(u/v), \quad \mbox{is completely monotone  in the variable}\; w = v + 1/v \ge 2,$$
(it is easy to see that $H_u$ is always a function of $w$).   In \cite[Theorem 5.3.1]{B}, it is shown that the class $\hcm$ corresponds to pointwise limits of functions of the form
\begin{equation} \label{dens}
x\; \mapsto\; C \;x^{a - 1}\prod_{i = 1}^N (x + c_i)^{-t_i}, \quad  N\in \Nset,\; C,\; \;c_i,\;\;t_i>0,\; \;0<a<\sum_{i = 1}^N c_i.
\end{equation}
Property  \cite[iv) pp. 68]{B} asserts that
\begin{equation}\label{propel}
f\in \hcm \Longleftrightarrow x^{a}f(x^b) \in \hcm, \;\;\mbox{for all $\; a\in \Rset$, and $|b|\leq 1$}.
\end{equation}
By Bondessons's definition \cite[Definition 9.10]{SSV},   and by \cite[Theorem 6.1.1]{SSV}, we have
\begin{equation}\label{gcgc}
X\sim \GGC \Longleftrightarrow \phi_X \in \TB \Longleftrightarrow e^{-\phi_X}\in \hcm.
\end{equation}
The class of $\hcm$ functions being stable under multiplicative convolution (\cite[property vii) p.68]{B}), then the Laplace transform of an $\hcm$ function is itself $\hcm$. If the p.d.f. of a continuous and positive r.v. $X$ is $\hcm$, we denote $X\sim\HCM$.  Clearly, we have the classification
$$\HCM \subset \GGC \subset  \BO, \;\;  \GGC \subset  \SD  \subset \ID,$$
where $\SD$ denotes the well-known class of {\it self-decomposable distributions}, see \cite{steutel} for this class.
Note that all previous classes of distributions are stable by additive convolution.
The importance of the  classes $\HCM$ and $\GGC $ is due to their specific stability properties, one can find in \cite{B1}: if $X$ and $Y$ are  independent and $ X, \; Y\sim \GGC$ (respectively $\; \HCM $), then
\begin{equation}\label{xx}
 X+Y,\;\; X \times Y\sim \GGC  \;\;\mbox{(resp. $\;X^q, \;q>1,\;\;\; X \times Y, \;\; X / Y  \sim \HCM$).}
\end{equation}
Note that the above stability properties are not shared by general distribution of the half real line. Besides, by \eqref{bb}, the $\ID$ property reads only on the level of the cumulant function definition, whereas the $\HCM$ one reads on the level of the probability density function, cf. the examples of next subsection.
\subsection{Gamma and positive stable distributions}
In what follows, $\ga_t$ denotes a standard gamma distributed random variable with shape with parameter $t>0$, i.e. with Laplace transform, Mellin transform and p.d.f.
\begin{equation}
\er[e^{-\lambda \ga_t}]=\Bigl(\frac{1}{1+\lambda }\Bigr)^t,\;\;
\lambda \geq 0,\quad \er[\ga_t^{s}]=\frac{\Gamma(t+s)}{\Gamma(t)} \;\;s>-t, \quad \mbox{and}\quad
f_{\ga_t}(x)=\frac{x^{t-1}}{\Ga(t)} e^{-x},\;\;x>0.\label{gamma}
\end{equation}
The last p.d.f. form and \eqref{dens} ensure that all gamma distributions are $\HCM$. Another example  of a $\GGC$ distribution is given by the {\it positive stable} r.v. $\sal, \;\alpha \in (0,1)$, which is associated to the Thorin Bernstein function
\begin{equation}\label{stal}
\phi_{\sal}(\lambda)=-\log \er[e^{-\lambda \sal}]= -\log\int_0^\infty e^{-\lambda x }\;f_{\sal} (x) \;dx =\lambda^\alpha=\frac{\alpha\sin(\alpha\pi)}{\pi}\int_0^\infty\log\left(1+\frac\lambda x\right) \; \frac{dx}{x^{1-\alpha}},\quad \lambda \geq 0.
\end{equation}
Note that the p.d.f.  $f_{\sal}$ is not explicit except for $\alpha =1/2$;
\begin{equation} \label{half}
\mathbb{S}_{1/2}\simdis \frac{1}{4 \;\ga_{1/2}}  \quad \mbox{and}  \quad f_{\mathbb{S}_{1/2}}(x)\; =\; \frac{1}{2\sqrt{\pi\, x^3}}e^{-1/4x},\;\; x>0,
\end{equation}
nevertheless, the Mellin transform of $\sal$ is explicitly given by
\begin{equation}\label{moms}
\er\left[ \sal^{-s} \right]=\frac{\Ga(1+\frac{s}{\alpha})}{\Ga(1+s)}, \quad s>-\alpha,
\end{equation}
which gives the convergence in distribution
\begin{equation}\label{t01}
\frac{1}{(\sal)^\alpha} \limdis 1, \;\; \mbox{as} \;\alpha \to 1,\quad \mbox{and} \quad  \frac{1}{(\sal)^\alpha} \limdis \exs, \;\; \mbox{as} \; \alpha \to 0.
\end{equation}
It is then natural to adopt the conventions $\mathbb{S}_0 =0,\;\;\;\mathbb{S}_1$. Related literature for the stable distribution is referenced on Nolan's website \cite{Nolan}. The monograph of Zolotarev \cite{zolo} is one of the major references for general stable distributions. By representation \eqref{stal}, and by the main result of \cite{BS}, we know that
\begin{equation}\label{stal121}
\sal \sim \GGC,\;\;\mbox{for all $\alpha\in (0,1),\quad$ and}\quad  \sal \sim \HCM\Longleftrightarrow\alpha \leq 1/2.
\end{equation}
Additionally, in \cite{JS}, it was shown  that if $t>0$, then
\begin{equation}\label{boto}
\sta^{-\frac{\al}{t}}\sim \ID \Longleftrightarrow  t\leq 1-\alpha.
\end{equation}
Using \eqref{t01}, and taking  $\sal'$  an independent version of $\sal$, let us define the  quotient, \begin{equation}\label{tuma}
\tal :=\left(\frac{\sal}{\sal '}\right)^\alpha,\;\;\alpha\in (0,1), \quad \mathbb{T}_1=1,  \quad\mbox{and}\quad   \mathbb{T}_0\simdis\frac{\exs}{\exs'},
\end{equation}
where $\exs,\,\exs '$ are i.i.d. and standard exponentially distributed. The  p.d.f. of $\tal$  is explicitly given by the generalized Cauchy form
\begin{equation}\label{stst}
 f_{\tal} (x) = \frac{\sin (\pi \alpha)}{\pi \alpha} \frac{1}{1+ 2 \cos(\pi \alpha)x +x^2},\quad x>0,
\end{equation}
see \eqref{momt} below for instance. Bosch \cite{bosch} completed \eqref{boto} by showing
\begin{equation}\label{bobo}
\alpha \in (0,1/2]\;\;\mbox{and} \;\;|t|\leq 1-\alpha\Longleftrightarrow  \tal^{1/t}=\left(\frac{\sal'}{\sal}\right)^{\alpha/t}\sim \HCM.
\end{equation}
With the above information,  Bosch and Simon found it natural to raise the following open question in \cite{BS}:
\begin{equation}\label{conju}
\mbox{\it Is it true that $\mathbb{S}_\alpha^{\alpha/t} \sim \HCM \;\;(\Longleftrightarrow\;\mathbb{S}_\alpha^{-\alpha/t} \sim \HCM)$ if, and only if,  $\alpha\leq 1/2$ and $|t|\leq 1-\alpha $?}
\end{equation}
This open question extends Bondesson's conjecture \cite{B}, which was stated in 1977. This conjecture, which  is identical to \eqref{conju} with $t=\alpha$, was investigated in \cite{JS} and solved in \cite{BS}.
\begin{figure}[h]
\begin{center}
\sidecaption
\includegraphics[scale=.8]{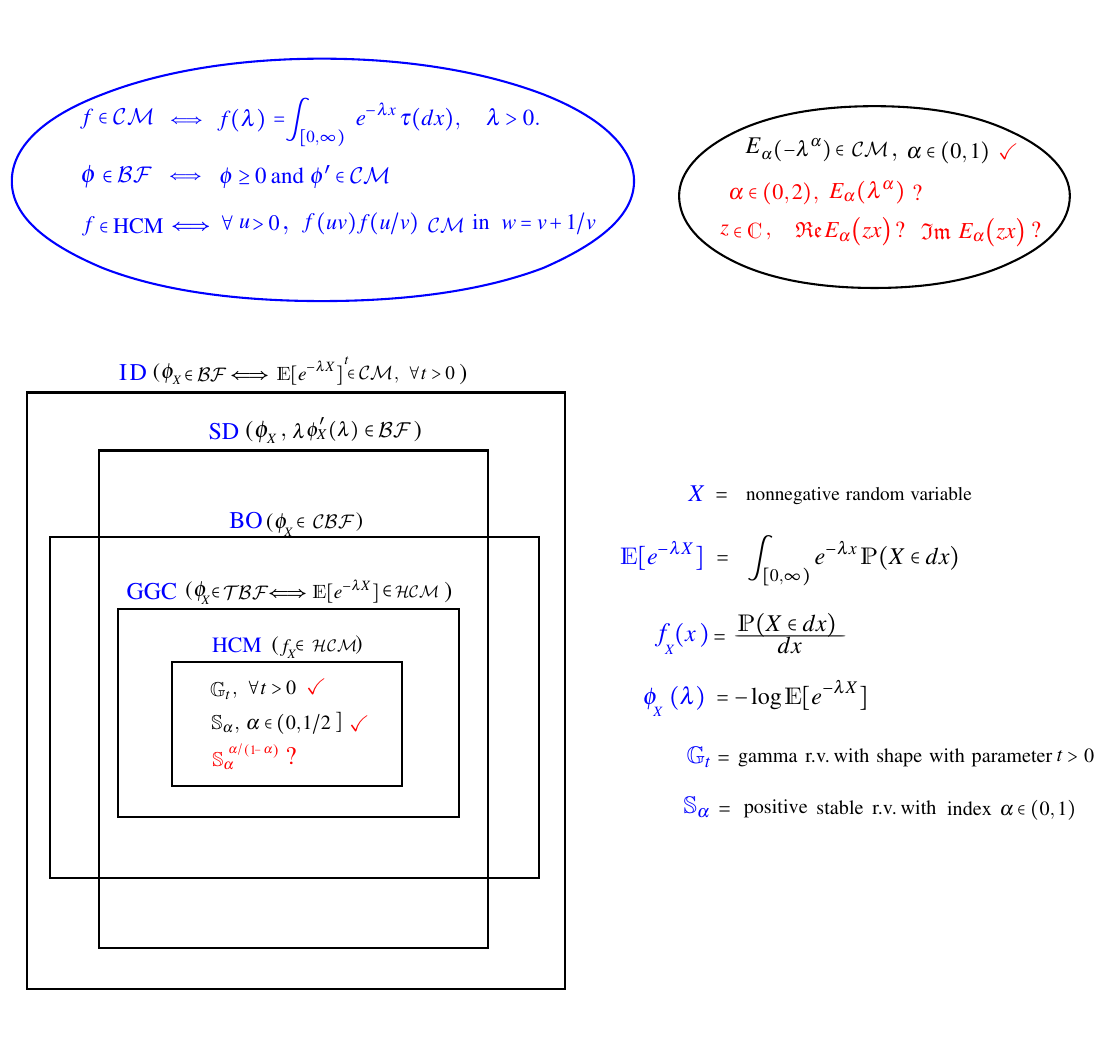}
\caption{Classes of infinitely divisible distributions, and motivations.}
\label{fig:Classes of ID}
\end{center}
\end{figure}
\subsection{First results on the Mittag-Leffler functions}

From \eqref{moms}, one gets a first link of the Mittag-Leffler function with the positive  stable  distribution:
\begin{equation}\label{mittag}
\er\Big[e^{z/(\sal)^{\alpha}}\Big]=E_\alpha (z), \quad \alpha \in (0,1), \; z\in \mathbb{C}.
\end{equation}
By  \eqref{stal}, it is clear that
\begin{equation}\label{mittagt}
\er[e^{-\lambda \tal^{1/\alpha}}]=\er[e^{-\lambda \sal'/\sal}]= \er[e^{-(\lambda/\sal)^\alpha}] =  E_\alpha (-\lambda^\alpha),\quad \lambda\geq 0\,.
\end{equation}
In \cite{JS}, if was shown that
\begin{equation*}\label{bibi}
0<\alpha \leq \frac{1}{2}\Longleftrightarrow \tal^{1/\alpha} \sim \HCM \Longrightarrow E_\alpha (-\lambda^\alpha) \in \hcm.
\end{equation*}
Note that the r.v.
\begin{equation}\label{tal}
\mathbb{R}_\alpha:=1+ \mathbb{T}_{|\alpha-1|} ^{1/\alpha}, \;\; \alpha \in [0,2],
\end{equation}
appeared in  \cite[(184) and (223)]{lry}, in a $\GGC$ context without noticing its $\GGC$ property: by \eqref{bobo}, we have
\begin{equation}\label{noti}
\frac{1}{2} \leq\al\leq 1 \Longleftrightarrow \mathbb{T}_{|\alpha-1|} ^{1/\alpha}=\mathbb{R}_\alpha-1\sim \HCM \Longrightarrow \mathbb{R}_\alpha \sim \GGC.
\end{equation}
Another monotonicity property of the Mittag-Leffler functions can be found in \cite{BJH}:
$$0<\alpha \leq \frac{1}{2}\Longleftrightarrow \lambda E_\alpha (-\lambda)\in \BF\Longleftrightarrow 1-   \,\Gamma(1-\alpha)\;\lambda\;E_\alpha(-\lambda)\in \CM,$$
and results in case $\alpha\in (1/2,1)$ are also shown there. Note  that the function $E_\alpha (-\lambda)$ is completely monotone, at  the contrary of $E_\alpha (\lambda)$.  In \cite[Theorem 1.2, (a) and (c)]{SIM3},  Simon showed that
\begin{equation}\label{simos}
\lambda \mapsto F_\alpha(\lambda):=\frac{e^\lambda - \alpha E_\alpha(\lambda^\alpha)}{1-\alpha} \in\CM, \;\;\mbox{for all} \;  \alpha \in(0,1)\cup(1,2].
 \end{equation}
The choice of the function $F_\alpha$  is intuitive because $ e^\lambda=E_1(\lambda) =\lim_{\alpha \to 1}\alpha E_\alpha(\lambda^\alpha)$.  \\

The main objective of this work is to exhibit new monotonicity  properties of the $\GGC$ and $\HCM$ type for the Mittag-Leffler functions. To give a taste of our results, we start by improving Simon's results \eqref{simos}.
\begin{teo} \label{mi1} For $\alpha\in (0,1)\cup(1,2]$, the following holds.
\begin{enumerate}[1)]
\item  With $\mathbb{T}_{|\alpha -1|}$  given by  \eqref{tuma}, and  $F_\alpha$ by \eqref{simos}, we have the Laplace representation:
\begin{equation}\label{Mittag1}
F_\alpha(\lambda):=
\er[e^{-\lambda {\mathbb{T}_{|\alpha -1|}}^{1/\alpha}}], \quad  \lambda \geq 0.
\end{equation}
\item  Assume $\alpha \in(0,1)$.
    \begin{enumerate}[(a)]
     \item The functions  $1- \alpha \;e^{-\lambda} E_\alpha (\lambda^\alpha)$ and $e^\lambda/E_\al(\lambda^\al)$ are completely monotone, and $e^{-\l} E_\alpha (\lambda^\alpha)-1$ is Bernstein.
    \item The function $F_\alpha$ is  Stieltjes.  Moreover, $F_\alpha \in \hcm$, if and only if, $1/2\leq \alpha <1$. In this case, $1- \alpha \;e^{-\lambda} E_\alpha (\lambda^\alpha)\in \hcm$  and  $-F_\alpha'/F_\alpha \in\ST_1$.
        \end{enumerate}
\item  Assume $\alpha \in(1,2]$.
    \begin{enumerate}[(a)]
    \item The function  $\alpha e^{-\l} E_\alpha (\lambda^\alpha)-1$ is completely monotone and $1- e^{-\l} E_\alpha (\lambda^\alpha)$ is Bernstein.
    \item The function $F_\alpha$  is not  $\hcm$  if $\alpha >4/3$.
    \item If $1<\alpha \leq 3/2$ and $0<\gamma \leq (2-\alpha)/\alpha$, then  $F_\alpha(\lambda^\gamma)$ is $\hcm$, and both functions $F_\alpha(\lambda^\gamma) $ and $-\lambda^{\gamma-1}\;F_\alpha'(\lambda^\gamma)/F_\alpha(\lambda^\gamma)$ are in $\ST_1$.
     \end{enumerate}
\end{enumerate}
\end{teo}

In the same direction as \eqref{simos}, Simon considered in  \cite{SIM1} and \cite{SIM3} the function
\begin{equation}\label{dda}
G_\al(x):=E_\al(x^\al)-\frac d{dx}\pab{E_\al(x^\al)}= \frac{1-\al}{\al}\big(F_\al'(x)-F_\al(x)\big), \;\;\;x>0.
\end{equation}
From \eqref{tal} and \eqref{Mittag1}, we deduce  the expressions
 \begin{equation}\label{-d}
G_\al(x)= \frac{\al-1}\al\;\er\left[\pa{1+\mathbb{T}_{|\alpha-1|}^{1/\alpha}} \;e^{-x \mathbb{T}_{|\alpha-1|}^{1/\al}} \right], \;\;\; \alpha\in (0,1)\cup (1,2].
\end{equation}
and
\begin{equation}\label{hh}
\frac{\al}{\al-1}\;e^{-x} G_\al(x)= \er\left[\mathbb{R}_\alpha \;e^{-x \mathbb{R}_\alpha} \right]\;\;\mbox{is the p.d.f. of the independent quotient $\ga_1/\mathbb{R}_\alpha$.}
\end{equation}
Using \eqref{moms}, we see that
\begin{equation}\label{d0}
\al\in(1,2]\Longrightarrow\er\left[\T_{\al-1}^{1/\al}\right]=\frac{1}{\alpha -1} \Longrightarrow G_\al(0+)=\frac{\al-1}\al\er\left[1+\T_{\al-1}^{1/\al}\right]=1= \int_0^\infty G_\al(x)\; dx,
\end{equation}
which means that  $G_\al, \; \al\in(1,2],$ is a p.d.f., and is also the Laplace transform of a probability distribution, see \eqref{abov} below.
In case $\alpha \in (0,1)$, we have $\er[\tumal ^{1/\alpha}]=\infty$ and the function $-G_\al$ is not a p.d.f. We complete our previous discussion with the following result.
\begin{teo}\label{Th:MittagLeffler2GGC} The r.v. $\mathbb{R}_\alpha$ and the function $G_\al$, given in \eqref{dda} and  \eqref{tal}, satisfy the following.
\begin{enumerate}[1)]
\item If $\al\in(1,2]$, then the p.d.f. $G_\al$ is the one of a $\GGC$ distribution but  not a $\HCM$.
\item If $\al\in(0,1)$, then the function $-G_\al$ is a widened $\GGC$, if and only if,  $\al\ge1/2$.
    In this case, $\ga_1/\mathbb{R}_\alpha\sim \GGC$.
\end{enumerate}
\end{teo}
\medskip

In \cite{SIM2}, Simon was interested in the first passage time of the normalized spectrally positive  $\alpha$-stable L\'evy process ${(S_t^{(\alpha)})}_{t\geq 0}, \;\alpha\in (1,2]$, viz.
$$\er[e^{-\lambda S_t^{(\alpha)}}]=e^{t\,\lambda^\alpha},\;\; \;\lambda, \;t\geq 0,$$
and
\begin{equation}\label{tau1}
\widehat{\tau}_x =\inf\{t >0\;;\; S_t^{(\alpha)}= -x\} \;\;\;\mbox{and}\;\;\; \tau_x =\inf\{t >0\;;\; S_t^{(\alpha)}= x\}, \;x>0.
\end{equation}
The self-similarity property of ${(S_t^{(\alpha)})}_{t\geq 0}$ with index $1=1/\alpha$, entails
$$S_t^{(\alpha)}\simdis t^{1/\alpha}\;S_1^{(\alpha)}, \quad (\widehat{\tau}_x, \tau_x)\simdis x^\alpha (\widehat{\tau}_1,\tau_1).$$
Additionally, it is well known that
$$\widehat{\tau}_1 \simdis \inf\{t >0\;;\; S_t^{(\alpha)}>1\}\simdis \mathbb{S}_{1/\alpha},$$
and this explains our focus on $\tau_1$. To formalize our next results, we need to recall the biasing procedure for the distribution of a non-negative r.v. $Z$. For $u \in \Rset$ such that $\er[Z^u]<\infty$, we denote by $Z^{[u]}$ a version of the length-biased distribution of order $u$  of $Z$, viz.
\begin{equation} \label{biais}
\pr(Z^{[u]}\in dx) \simdis   \frac{x^{u}}{\er[Z^u]} \;  \pr(Z\in dx).
\end{equation}

 We complete Simon's main result in \cite{SIM2} by expliciting (in \eqref{tau2} below) his factorization of $\tau_1$ obtained in his main theorem,  and by the use of the function
\begin{equation}\label{simohs}
H_\alpha(\lambda):=(1-\alpha)F_\alpha'(\lambda^{1/\alpha})=e^{\lambda^\frac{1}{\alpha}} -\alpha^2 \;\lambda^{1-\frac{1}{\alpha}}\; E_\alpha'(\lambda), \;\;\;\alpha \in (1,2].
\end{equation}
\begin{cor} Let  $\alpha \in (1,2]$. The r.v. $\tau_1$ given by \eqref{tau1}, and $H_\alpha$ by \eqref{simohs}, satisfy the following.
\begin{enumerate}[1)]
\item We have the Laplace representation:
\begin{equation}\label{galga}
\er[e^{-\lambda \tau_1}] = H_\alpha(\lambda), \;\; \lambda \geq 0.
\end{equation}
\item  With the convention \eqref{tuma}, we have  the independent factorizations in law:
\begin{equation}\label{tau2}
\tau_1 \simdis \mathbb{S}_{1/\alpha} \;(\talm)^{\;[1/\alpha]}\quad \mbox{and}\quad
\left(\frac{\ga_1}{\tau_1}\right)^{1/\alpha} \simdis \ga_1\; (\talm^{1/\alpha})^{\;[1]},
\end{equation}
where the distribution of  $(\talm)^{\;[1/\alpha]}$ and $(\talm^{1/\alpha})^{\;[1]}$ are given by the biasing procedure \eqref{biais}.
\item $(\ga_1/\tau_1)^{1/(2-\alpha)}\sim \HCM$.
\item If  $\alpha \leq 3/2$, then $\tau_1\sim \GGC$, hence $H_\alpha  \in \hcm.$
\end{enumerate}
\label{to1} \end{cor}

In addition to the above results, which  are a  continuation of Simon's ones  \cite{SIM1,SIM2,SIM3}, we obtain the following ones that we divide, for sake of coherence, into Sections \ref{addi} and \ref{Stochastic}:
\begin{itemize}
\item Proposition \ref{hcmlog} is a key to illustrate several other links between the r.v.'s $\tal$ \eqref{tuma}, and the Mittag-Leffler functions;
\item A consequence of Proposition \ref{hcmlog} is Corollary \ref{hcmcor}, that builds a new class of $\hcm$ functions, similar to the shapes \eqref{dens}  and \eqref{stst};
\item In Corollary \ref{crt1} (and also in Theorem \ref{sin}), we obtain that, for any complex number $z$ in the first quadrant, and for $0<t\leq 1-\arg z$, the real and imaginary part of  $\lambda\mapsto E_t \big(-z \lambda^t\big)$ also enjoy  $\CM$, $\GGC$ and $\HCM$ properties;
\item  With the help of Proposition \ref{hcmlog},  we characterize in Theorem \ref{sin} the distributions of a peculiar family of distributions, that might be helpful in solving the open  question \eqref{conju};
\item  We conclude with Corollary \ref{corg}, which provides more information that \eqref{bobo}.
\end{itemize}
To the best of our knowledge, the results of Sections \ref{addi} and \ref{Stochastic} have not been addressed in the literature in a form even close to ours, and they constitute a contribution to infinite divisibility, especially in the Bondesson's $\HCM$  and Mittag-Leffler contexts.

In Section \ref{setting}, we  give  comments and perspectives on some of the above results and  also settle a prerequisite on the positive stable distributions on the Bernstein functions. The proofs are postponed to Section \ref{theproofs}.
\section{A new class of $\HCM$ distributions and $\CM$ property for $\lambda \mapsto\Rel \big(E_t(-z\lambda^t)\big), \;z\in \Cset$.}\label{addi}

\subsection{The biasing  and the gamma-mixture procedure.}
Clearly,
\begin{equation} \label{biaiso}
\er[(Z^{[u]})^\lambda]= \frac{\er[Z^{\lambda+u}]}{\er[Z^{u}]} \quad\mbox{and}\quad (Z^v)^{[u]}\simdis (Z^{[uv]})^v, \;\mbox{if}\;\;\er[Z^{uv}]<\infty.
\end{equation}
A  useful result is the stochastic interpretation of property  \eqref{propel}:
\begin{equation}\label{prov}
Z\sim \HCM \Longleftrightarrow (Z^q)^{[u]}\simdis (Z ^{[qu]})^q\sim \HCM, \;\;\mbox{for  $|q|\geq 1$ and all $u$ s.t.}\;\er[Z^{uq}]<\infty.
\end{equation}
A random variable $\Be_{a,b},\,a,\,b>0$,  has the beta distribution  if it has the p.d.f.
$$\frac{\Gamma(a+b)}{\Gamma(a)\,\Gamma(b)} \,x^{a-1}(1-x)^{b-1}, \quad 0<x<1.$$
Observe that $\ga_t^{[u]}\simdis \ga_{t+u}$ and $\Be_{a,b}^{[v]}\simdis \Be_{a+v,b}$, if $u>-t$ and $v>-a$. By the {\it beta-gamma algebra}, we know that  if $\Be_{a,b}$ and $\ga_{a+b}$ are independent, then
\begin{equation}
(\Be_{a,b}, \ga_{a+b})  \simdis  \left(\frac{\ga_u}{\ga_a+\ga_b}, \ga_a+\ga_b\right),\quad \mbox{where $\ga_a$ and $\ga_b$ are independent}.
\label{alge} \end{equation}
The following fact will be  used  in the sequel to exhibit several factorizations in law. A positive r.v. $X$ has p.d.f. of the form
\begin{equation} \label{factuel}
f_X(x) = x^{t-1} g(x),\;\; t>0 \;\; \mbox{and}\; g \in \CM,
\end{equation}
if, and only if, the distribution of $X$ is a {\it gamma-mixture} of order $t$ (shortly $X$ is  a $\ga_t$-{\it mixture}). The latter is equivalent to the existence of a  positive r.v. $Y$ such that we have the independent factorization
\begin{equation} \label{fact}
X\simdis \frac{\ga_t}{Y},
\end{equation}
and in \eqref{fact}, necessarily
$$g(x)= \frac{1}{\Ga(t)}\;\er\left[Y^t\; e^{-xY}\right]\; \left(=\frac{\er[Y^t]}{\Ga(t)} \; \er\big[\; e^{-xY^{[t]}}\big],\;\;\mbox{if $\er[Y^t]<\infty$}\right).$$
An important result due to Kristiansen \cite{Kri} asserts that
\begin{equation}\label{infig}
X\; \mbox{\it has a $\ga_2$-mixture distribution} \Longrightarrow X\sim \ID,
\end{equation}
and the beta-gamma algebra ensures that the same holds if the distribution of $X$ is a $\ga_t$-mixture, $t<2$.
\subsection{A generalization of property \eqref{bobo}.}\label{first}
In this work, the functions, defined for $c \geq 0$ and $t>0$, by
\begin{equation}\label{psia}
\phi_{c,t}(\lambda):=\log\left(1 + 2 \;c \;\lambda^t+  \lambda^{2t}\right), \quad \lambda \geq 0,
\end{equation}
turn out to be crucial, and their properties will unblock several of our problems. Notice that the necessity of the condition $c\geq 0$ is essential since we need $\phi_{c,t}$   to be a Bernstein function  in the sequel. Notice  that $\phi_{c,t}$ is the composition $\phi_{c,t}(\lambda)=\varphi(\psi_{c,t}(\lambda))$, where $\varphi(\lambda):=\log(1+\lambda)\in \TB$, whereas
$$\psi_{c,t}(\lambda):=  2 \;c \;\lambda^t+  \lambda^{2t}\in \TB \;\;\mbox{if, and only if, $t\leq 1/2$}.$$
On the other hand, the class $\CB$ being stable by composition,   we see that
$$\lambda \mapsto \frac{\lambda^{2t}}{ 2\;c\;\lambda^{t}+\lambda^{2t}}\;\; \mbox{is a complete Bernstein function  if $t\leq 1/2$}.$$
By \cite[Theorem 7.3]{SSV}, we deduce that the logarithmic derivative
$$\frac{\psi_{c,t}'(\lambda)}{\psi_{c,t}(\lambda)}= \frac{t}{\lambda}\left(1+ \frac{\lambda^{2t}}{ 2\;c\;\lambda^{t}+\lambda^{2t}}\right) \;\;\mbox{is a Stieltjes function  if}\;\; t\leq 1/2,$$
i.e. $\psi_{c,t}'$ has the representation by \eqref{repsti}. By \eqref{compotb}  below, we easily deduce that $\phi_{c,t} \in \TB$,  for any $c\geq 0$ and $t\leq 1/2$. A larger range for $t$ is provided by the following result.
\begin{prop}  Let $t > 0$ and $c\geq 0$.  Then the following holds.
\begin{enumerate}[1)]
\item The function
\begin{equation}\label{logt}
\phi_{c,t}(\lambda):=\log (1+2\; c \lambda^t + \lambda^{2t}), \quad \lambda\geq 0,
\end{equation}
is Thorin Bernstein if, and only if, $t\leq 1$ and $c+\cos (\pi t)\geq 0$.
\item Let $\alpha \in [0,1]$ and $c=\cos(\pi \alpha)$. Then the following assertions are  equivalent.
      \begin{enumerate}[(i)]
        \item $\alpha \leq 1/2$ and $t\leq 1-\alpha$;
        \item The function $(1+2 \cos(\pi \alpha)\lambda^t + \lambda^{2t})^{-1}$ is completely monotone;
        \item The function $\phi_{c,t}$ is Bernstein;
        \item The function $\phi_{c,t}$ is Thorin Bernstein.
        \end{enumerate}
\end{enumerate}
\label{hcmlog}\end{prop}
\medskip

A consequence  of Proposition \ref{hcmlog} is the following result which gives more information on
the function $\phi_{c,t}(\lambda)$, and provides a monotonicity property for the Mittag-Leffler function in case $c\in [0,1]$.
\begin{cor}  Let $c\in [0,1], \;\jc := c+\ii\sqrt{1-c^2}.$ For $t\in(0,1)$, let $E_t$ be the Mittag-Leffler function and $\sat$ be a positive stable random variable. Then, the following holds.
\begin{enumerate}[1)]
\item The series
\begin{equation}\label{surp}
\mathfrak{C}_{c,t}(x):= \sum_{n=0}^{\infty} (-1)^n \;\frac{\cos\big(n \arccos(c)\, \big)}{\Ga\big(nt+1\big)}\; \;x^{nt},\;\;\;  x\geq 0,
\end{equation}
is represented by
$$\mathfrak{C}_{c,t}(x) = \Rel \left(E_t \big(- \jc\;x^t\big)\right)=\er\left[\chi_c\left(\Big(\frac{x}{\sat}\Big)^t\right)\right],  \;\;\;  \mbox{where}   \;\;\; \chi_c(x):= e^{-cx}\cos(\sqrt{1-c^2}\, x).$$
\item The following assertions are equivalent.
   \begin{enumerate}[(i)]
   \item $0\leq t \leq 1-\arccos(c)$;
   \item the function $\mathfrak{C}_{c,t}$ is completely monotone;
   \item the function $1-\mathfrak{C}_{c,t}$ is Bernstein;
   \item the function $\phi_{c,t}(\lambda)=\log(1 + 2 \;c \;\lambda^t+  \lambda^{2t})$ is Thorin Bernstein.
    \end{enumerate}
    \smallskip
\item Under any of the conditions in 2),  the functions $\phi_{c,t}$ and $\mathfrak{C}_{c,t}$ are represented by
     \begin{equation}\label{pip}
    \phi_{c,t}(\lambda)=2 t \int_{0}^{\infty}(1-e^{-\lambda x})\frac{\mathfrak{C}_{c,t}(x)}{x}dx, \quad \mathfrak{C}_{c,t}(y)=\er[e^{-y \;\mathbb{E}_{c,t}}], \quad \lambda, \;y \geq 0,
    \end{equation}
     where $\mathbb{E}_{c,t}$ is a positive r.v. whose distribution gives no mass in 0. The Thorin measure $U_{c,t}$ of $\phi_{c,t}$, obtained by representation \eqref{repphi3tt}, is  $U_{c,t}(du)=2 \;t \; \,\pr(\mathbb{E}_{c,t}\in du)$.

\item Let $z$ be a complex number such that $\Re(z),  \;Im(z)\geq 0$. Then, $\lambda \mapsto \Rel \left(E_t \big(- z\;\lambda^t\big)\right) \in \CM$, if, and only if, $t\leq 1-\arg(z)$.
\end{enumerate}

\label{crt1} \end{cor}
\begin{rmk}  Using \eqref{mittag}, \eqref{mittagt}  and \eqref{surp}, we obtain
$$\mathfrak{C}_{1,t}(x)= \er\big[e^{- (x/\sat)^t}\big]= \er\big[e^{-x\;\tat^{1/t}}\big] = E_t(-x^t) \quad \mbox{and}\quad\mathfrak{C}_{0,t}(x)= \er\left[\cos\big(x^t/\sat^t\big)\right]= \Rel\, E_t (\ii\;x^t),$$
then the completely monotonicity of   $\mathfrak{C}_{1,t} $ is not a surprise. On the other hand, Proposition \ref{hcmlog} gives
\begin{eqnarray*}
0<t\leq 1& \Longleftrightarrow&
\phi_{1,t}(\lambda)=2\;\log(1+\lambda^{t})= 2\;t \int_0^\infty  \!\!(1- e^{-\lambda x}) \;\frac{ E_t(-x^t)}{x}\;dx \in \TB \Longleftrightarrow E_t(-x^t)  \in \CM ,\\
0<t\leq 1/2&\Longleftrightarrow & \phi_{0,t}(\lambda)= \log(1+\lambda^{2t})= 2\;t \int_0^\infty \!\! (1- e^{-\lambda x}) \;\frac{\Rel\big( E_t (\ii\;x^t)\big)}{x}\;dx \in \TB\Longleftrightarrow  \Rel\big( E_t (\ii\;x^t)\big)\in \CM.
\end{eqnarray*}
The trivial relation
$$\Rel\, E_a (\ii\;x)= \sum_{n \geq 0} (-1)^n \frac{x^{2n}}{\Gamma (2n a +1)} = E_{2\,a}(-x^2), \quad x\geq 0, \;\;a>0,$$
complies with  the equality $\phi_{1,t}(\lambda)=2\;\phi_{0,t/2}, \;t \leq 1$.
\end{rmk}


We are now able to introduce a new class of $\hcm$ p.d.f.'s that are reminiscent of representation \eqref{dens}.

\medskip
\begin{cor} Probability density functions of the form
\begin{equation}\label{newp}
f_{c,a,t}(x)=C\; \frac{x^{a-1}}{1+2\; c \;x^t + x^{2t}}, \quad C, \; t,\;a, \;x>0,
\end{equation}
are $\hcm$, if, and only if,  $t\leq 1$, $c\geq 0$ and $c+\cos (\pi t)\geq 0$.
Pointwise limits of functions of the form
$$ C\; x^{a}\prod_{1}^{N}\big(1+2\; c_i \;x^{t_i} + x^{2t_i }\big)^{-\gamma_i}, \quad N\in \Nset, \;C,\; \gamma_i>0,\;\; a\in \Rset,\;\;  0< t_i\leq 1,  \;\;c_i\geq 0,\;\;\mbox{and}\;  \;c_i+\cos (\pi t_i)\geq 0,$$
are $\hcm$.
\label{hcmcor}\end{cor}
\
\medskip
Thanks to \eqref{stst},  the p.d.f. of $\;\tal^{1/t}= (\sal'/\sal)^{\alpha/t}, \; t>0$, is explicitly given  by
\begin{equation}\label{ftal}
f_{\tal^{1/t}}(x) = t \,x^{t-1}\, f_{\tal}(x^t)=\frac{t \sin(\pi \alpha)}{\pi \alpha}\frac{x^{t-1}}{1+2 \;\cos(\pi \alpha) \;x^t + x^{2t}}, \quad x>0.
\end{equation}
Therefore, the main result of Bosch \cite[Theorem 1.2]{bosch}  is a particular case of Corollary \ref{hcmcor} with $\;c=\cos(\pi \alpha)$ and $a=t$ there. Indeed, the condition $0<t\leq 1$, $c\geq 0$ and $c+\cos (\pi t)\geq 0$ is equivalent to
\begin{equation}\label{fou}
0\leq \alpha \leq \frac{1}{2}, \;\;\;  0\leq t \leq 1-\alpha,
\end{equation}
and then \eqref{bobo} holds true.

\section{Stochastic interpretation of the p.d.f.'s in \eqref{newp} and $\CM$ property for $\lambda \mapsto\Ima \big(E_t'(-z\lambda^t)\big), \;z\in \Cset$.}\label{Stochastic}
Proposition \ref{hcmlog} enables us to introduce the positive r.v. $\mathbb{X}_{c,t} \sim \GGC$  associated with the $\TB$-function $\phi_{c,t}$ in  \eqref{logt}, i.e.
\begin{equation}\label{lala}
 \er[e^{-\lambda \mathbb{X}_{c,t}}]=\frac{1}{1+2c\lambda^t+\lambda^{2t}}=e^{-\phi_{c,t}(\lambda)}, \quad 0<t\leq 1, \;\;c\geq 0 \;\;\mbox{and}\; \;c+\cos(\pi t)\geq 0.
\end{equation}
With the convention $\mathbb{S}_1=1$, the random variable  $\mathbb{X}_{1,t}$ (respectively for $\mathbb{X}_{0,t}$) is well defined for $t\leq 1$  (respectively $0<t\leq 1/2$). These r.v.'s
enjoy  a simple independent factorization in law:
\begin{equation}\label{x1}
\er[e^{-\lambda \mathbb{X}_{1,t}}]=\frac{1}{(1+\lambda^t)^2}=
\er[e^{-\lambda \mathbb{S}_t\;\ga_2^{1/t}}] \Longleftrightarrow\mathbb{X}_{1,t}\simdis \mathbb{S}_t \;\ga_2^{1/t},
\end{equation}
and
\begin{equation}\label{x0}
\er[e^{-\lambda \mathbb{X}_{0,t}}]=\frac{1}{1+\lambda^{2t}}=
\er[e^{-\lambda \mathbb{S}_{2t}\;\ga_2^{1/t}}] \Longleftrightarrow\mathbb{X}_{0,t}\simdis \mathbb{S}_{2t} \;\ga_1^{1/t}.
\end{equation}
For $0<t<s\leq 1$ and $c$  as in Proposition \ref{hcmlog}, we have this algebra: if  $\mathbb{S}_{t/s}$ is independent of $\mathbb{X}_{c,s}$, then
\begin{equation}\label{dup}
\er[e^{-\lambda \mathbb{X}_{c,t}}]= \er[e^{-\lambda^{t/s} \mathbb{X}_{c,s}}]= \er[e^{-\lambda \mathbb{S}_{t/s} (\mathbb{X}_{c,s})^{s/t}}]  \Longrightarrow \mathbb{X}_{c,t}\simdis \mathbb{S}_{t/s} \,(\mathbb{X}_{c,s})^{s/t}.
\end{equation}
The latter is a reminiscent of the {\it subordination relation} for stable distributions: if  $\mathbb{S}_{\alpha/\gamma}, \;0<\alpha<\gamma<1,$ is a stable r.v. independent of $\mathbb{S}_\gamma$, then
\begin{equation}\label{subo}
\er[e^{-\lambda \sal}]= e^{-  \lambda^\alpha}= e^{-  (\lambda^{\alpha/\gamma})^\alpha}=
\er[e^{-\lambda^{\alpha/\gamma} \sal}]= \er[e^{-\lambda \mathbb{S}_{\alpha/\gamma} (\mathbb{S}_\gamma)^{\gamma/\alpha}}] \Longrightarrow \sal \simdis \mathbb{S}_{\alpha/\gamma}\; (\mathbb{S}_\gamma)^{\gamma/\alpha}.
\end{equation}
The ordinary generating function for the Chebychev's polynomial of the second kind  is
$$\sum_{n\ge0} U_n(x)\,u^n\,=\,\frac1{1-2\;x\;u+u^2}, \quad |x|, \, |u|<1.$$
When $c\in [0,1]$, me may take
\begin{equation}\label{fou}
c=\cos (\pi \alpha),\;\; 0\leq \alpha \leq \frac{1}{2}, \;\;\;  0\leq t \leq 1-\alpha,
\end{equation}
and get the Laplace transform
\begin{equation}\label{yct}
\er[e^{-\lambda \mathbb{X}_{c,t}}]=\frac{1}{1+2\cos(\pi \alpha)\lambda^t+\lambda^{2t}} = \sum_{n\geq 0}U_n(-\cos \pi \alpha)\; \lambda^{nt}= \sum_{n\geq 0}(-1)^n \frac{\sin \pi (n+1) \alpha}{\sin \pi \alpha}\; \lambda^{nt},  \qquad |\lambda|<1.
\end{equation}
Using the convention \eqref{fou} and  the form \eqref{ftal}, we get
$$f_{\tal^{1/t}}(x)=t\frac{\sin(\pi \alpha)}{\pi \alpha}  x^{t-1}\;\er[e^{-x \mathbb{X}_{c,t}}], \;\;x>0.$$
Integrating the latter on $(0,\infty)$, we obtain
\begin{equation}\label{xit}
\er[(\mathbb{X}_{c,t})^{-t}]= \frac{\pi \alpha} {\sin(\pi \alpha)\Gamma(t+1)}.
\end{equation}
By \eqref{biais} and \eqref{fact}, we deduce
\begin{eqnarray*}
f_{\tal^{1/t}}(x)&=&   \frac{x^{t-1}}{\Gamma(t)}\;\er\left[ \left(\mathbb{X}_{c,t}^{[-t]}\right)^t e^{-x \mathbb{X}_{c,t}^{[-t]}}\right], \;\;x>0.
\end{eqnarray*}
Hence,
\begin{equation}\label{xatc}
\alpha\in (0,1/2], \; t\in (0,1-\alpha] ,\;\;\; \mbox{and}\;\;\;c=\cos(\pi \alpha) \Longrightarrow
\tal^{1/t}\simdis \frac{\ga_t}{(\mathbb{X}_{c,t})^{[-t]}}.
\end{equation}
Then taking the  Mellin transform in both sides of the latter, we get
$$\er[(\mathbb{X}_{c,t}^{[-t]}) ^{-\lambda}] =\frac{\er[(\mathbb{X}_{c,t}) ^{-(\lambda+t)}]}{\er[\mathbb{X}_{c,t}^{-t}]}=
\frac{\Gamma(t)\;\Gamma(1-\frac{\lambda}{t})\Gamma(1+\frac{\lambda}{t})} {\Gamma(t+\lambda)\;\Gamma(1-\frac{\alpha\lambda}{t})\;)\;\Gamma(1+\frac{\alpha\lambda}{t})} ,\quad  |\lambda|<t,$$
or equivalently
\begin{eqnarray}
\er[\mathbb{X}_{c,t}^{-x}]&=&\frac{\pi \alpha}{\sin(\pi \alpha)} \frac{\Gamma\left(2- \frac{x}{t}\right)\; \Gamma\left(1+\frac{x}{t}\right)}{\Gamma(1+x)\; \Gamma\left(1- \alpha(1-\frac{x}{t})\right)\; \Gamma\left(1+ \alpha(1-\frac{x}{t})\right)} \nonumber\\
&=&\left\{
\begin{array}{ll}
\displaystyle \frac{\pi \alpha} {\sin(\pi \alpha)\Gamma(t+1)},\quad & {\rm if }\;\; x =t\\
&\\
\displaystyle \frac{\Gamma(1-\frac{x}{t})\Gamma(1+\frac{x}{t})} {\Gamma(1+x)} \frac{\sin\big(\pi \alpha(1-\frac{x}{t})\big)}{\sin\big(\pi \alpha)},\quad & {\rm if }\; \; x\in [0,t)\cup(t,2t).
\end{array} \right. \label{carc}
\end{eqnarray}

Motivated by the  duplication formula \eqref{dup}, the link \eqref{yct} with Chebychev's polynomials, and \eqref{carc},  we explicit the distribution of the r.v.'s $\mathbb{X}_{c,t}, \; c\in(0,1)$, in the following result, the case $c=1,0,$ is described in \eqref{x1} and \eqref{x0}.
\begin{teo} Let $\alpha \in (0,1/2],  \;t\in  (0,1-\alpha]$,  $c=\cos(\pi \al), \; j_c=e^{\ii \pi \al}$ and $\jcb=e^{-\ii \pi \al}$ . The $\GGC$ distribution of  r.v. $\mathbb{X}_{c,t}$ is described by the following.
\begin{enumerate}[1)]
\item The c.d.f. and the p.d.f.  of the $\mathbb{X}_{c,t}$   are
$$\pr(\mathbb{X}_{c,t}\leq x)=\frac{1}{\sin(\pi \al)}\Ima \left(\frac{E_t(-\jc\, x^t)-1}{\jc}\right), \;\; x\geq 0,$$
\begin{equation}\label{well}
f_{_{\mathbb{X}_{c,t}}}(x)=\frac{t}{\sin(\pi \al)}\;x^{t-1}\;\Ima\big(-E_t'(-\jc\, x^t)\big)=\frac{t}{\sin(\pi \al)}\;x^{t-1}\; \Ima\left(\er\left[\frac{e^{-  \jcb(x/\mathbb{S}_t)^t} }{\mathbb{S}_t^t}\;
\right]\right), \;\;\; x>0.
\end{equation}
\item There exists a positive r.v. $\dac$, such that
$$\er[\dac^{2t}]=1, \quad \dac^{\;[t]} \simdis \frac{1}{\dac^{\;[t]}},$$
whose Mellin transform is
\begin{equation} \label{carcd}
\er\left[\left(\dac^{[t]}\right)^\lambda\right]=\frac{\Gamma(t)^2\;\Gamma(1-\frac{\lambda}{t})
\Gamma(1+\frac{\lambda}{t})} {\Gamma(t-\lambda)\;\Gamma(t+\lambda)\;\Gamma(1-\frac{\alpha\lambda}{t}) \;\Gamma(1+\frac{\alpha\lambda}{t})} ,\quad  |\lambda|<t,
\end{equation}
and such that
we have the independent factorization in law
\begin{equation}\label{ideng}
\mathbb{X}_{c,t} \simdis \frac{\ga_{2t}}{\dac}\sim \GGC,  \quad \mbox{and}\;\;\;(\mathbb{X}_{c,t})^{[-t]}\simdis \frac{\ga_t}{\dac^{\;[t]}} \sim \GGC,
\end{equation}
(recall the size biasing notation  \eqref{biais} for $Z^{[u]}$).  In particular, we have the Laplace transform representation:
\begin{equation}\label{wella}
\er[e^{-\lambda \dac^{[2t]}}]=\frac{t\;\Gamma(2t)}{\sin(\pi \al)}  \;\frac{\Ima\big(-E_t'(-\jc\, \lambda^t)\big)}{\lambda^t}=\frac{t\Gamma(2t)}{\sin(\pi\al)\;\lambda^t}\; \er\left[\frac{e^{-\cos(\pi\al)(\lambda/\mathbb{S}_t)^t}}{\mathbb{S}_t^t} \; \sin\big(\sin(\pi\al)\; (\lambda /\mathbb{S}_t)^t\big)\right],\;\;\lambda \geq 0.
\end{equation}
\item If $z$ is a complex number such that $\Re(z),  \;Im(z)\geq 0$ and $0<t\leq 1-\arg(z)$, then $\lambda \mapsto \Ima \big(E_t' (- z\;\lambda^t)\big)/\lambda^t \in \CM$.
 \end{enumerate}
\label{sin}\end{teo}

\medskip
The next result completes Theorem \ref{sin} and Bosch's characterization  \eqref{bobo}.
\begin{cor} For $\alpha\in (0,1)$ and $t>0$, we have the equivalences:
\begin{enumerate}[1)]
\item $\alpha\leq 1/2$ and  $t\leq 1-\alpha$;
\item  $x\mapsto (1+2 \;\cos(\pi \alpha) \;x^t + x^{2t})^{-1}\in \CM$;
\item  $x\mapsto (1+2 \;\cos(\pi \alpha) \;x^t + x^{2t})^{-1}\in \ST_{2t}$;
\item the distribution of $\tal^{1/t}$ is a $\ga_t$-mixture;
\item $\tal^{1/t} \;\sim\HCM$.
\end{enumerate}
With $c=\cos(\pi \alpha)$, we have the independent factorizations in law
\begin{equation}\label{idt}
\tal^{1/t}= \left(\frac{\sal}{\sal'}\right)^{\frac{\al}{t}} \simdis \frac{\ga_t}{(\mathbb{X}_{c,t})^{[-t]}}\simdis  \frac{\ga_t}{\ga_t'} \dac^{\;[t]}\sim \HCM.
\end{equation}
\label{corg}\end{cor}
\section{Comments and prerequisite for the proofs.} \label{setting}
\subsection{Comments on Theorems \ref{mi1}.}
Our results  in Theorems  \ref{mi1} merit  some comments.
\begin{enumerate}[(a)]
\item By \eqref{fact} and by \eqref{Mittag1}, it is clear that the function $-F_\alpha'$   is the p.d.f. of the independent quotient $\ga_1/\T_{|\al-1|}^{1/\al}$.
\item By \eqref{t01}, we have the convergence in law
\begin{equation}\label{latter}
\mathbb{T}_{|\alpha -1|}^{1/\alpha}\limdis \frac{\ga_1}{\ga_1'}, \quad\mbox{as}\;\;\al\to 1,
\end{equation}
where the independent ratio $\ga_1/\ga_1'$ has the Pareto with the $\hcm$ p.d.f. $(x+1)^{-2}, \;\;x>0$. Note that the latter can be obtained using \eqref{bobo} and \eqref{latter}, with $t=1-\alpha \to 0$, and the stability of the $\hcm$ class under pointwise limits. By \eqref{Mittag1}, the latter also reads
$$\lim_{\alpha \to 1}F_\alpha(\lambda)= \left(\frac{\partial \big(\al E_\alpha(\lambda^\alpha)\big)}{\partial\alpha}\right)_{|\alpha =1} = e^\lambda + \left(\frac{\partial\big( E_\alpha(\lambda^\alpha)\big)}{\partial\alpha}\right)_{|\alpha =1}=\int_{0}^{\infty} \frac{e^{-\lambda x}}{(x+1)^{2}}dx, \quad \lambda \geq 0.$$
\end{enumerate}
\subsection{Comments on Theorem \ref{Th:MittagLeffler2GGC}.}
The case $\alpha \in (1,2)$ in Theorem  \ref{Th:MittagLeffler2GGC} is specific.
\begin{enumerate}[a)]
\item  The complete monotonicity of the function $G_\alpha$ in \eqref{Form:DefinitionDal} is also a direct consequence of the one of $F_\alpha$ in \eqref{Mittag1}. Indeed,  since $F_\alpha\geq 0$, then
$$F_\alpha \in \CM \Longleftrightarrow -F_\alpha' \in \CM \Longrightarrow G_\alpha=F_\alpha-F_\alpha' \in \CM.$$

The functions
\begin{equation*}\label{obo}
-(\alpha-1)\; F_\alpha'(\lambda)= e^\lambda - \alpha \frac{d}{d\lambda} E_\alpha(\lambda^\alpha)= \er[e^{-\lambda^\alpha\; \tau_1}],\quad \alpha \in (1,2),
\end{equation*}
and  $F_\alpha, \; \alpha \in (0,1)\cup(1,2)$, were  shown to be completely monotone in \cite[(1.4)]{SIM2} and in \cite[Theorem 1.1]{SIM3}, respectively. The connection with the r.v. $\mathbb{T}_{1-\alpha},\;\alpha \in (0,1)$ was not noticed there, nor was the $\hcm$ property of the function $F_\alpha,\;\alpha\in (0,1)\cup(1,2)$.
\item  One has $-D_{1/2}(x)=(\pi\sqrt x)^{-1}\in \hcm$. One could ask if there exists some $\al\in(1/2,1)$ such that $-G_\al\in \hcm$. Observe that, contrary to the case $\al\in(1,2]$, the function
$$x\mapsto \frac{1+x}{x^{2\al}-2\cos(\pi\al)x^\al +1}, \quad x\geq 0,\;\; \al\in(1/2,1),$$
is completely monotone (one could even show that it is a Stieltjes function).\\
\item  In \cite[(3)]{SIM1} and \cite[(2.1)]{SIM3},   on can find the  following representations valid for $\alpha\in (0,1)\cup (1,2]$:
\begin{align}\label{Form:DefinitionDal}
G_\al(x)=-\frac{\sin(\pi\al)}\pi\int_0^\infty e^{-xu}\, \frac{u^{\al-1}(1+u)}{u^{2\al}-2\cos(\pi\al)u^\al+1}\, du,\;\;\; x>0,
\end{align}
and
\begin{equation}\label{Form:LaplaceDal}
\Lap_{G_\al}(\la)=\int_0^\infty e^{-\la x}G_\al(x)\,dx=\frac{\la^{\al-1}-1}{\la^\al-1},\quad \la\ge 0.
\end{equation}
In \eqref{d0}, we have seen that $G_\al(0+)=1$, in case $\alpha\in (1,2]$, consequently, the function
\begin{equation}\label{abov}
-\frac{\sin(\pi\al)}\pi\times\frac{u^{\al-1}(1+u)}{u^{2\al}-2\cos(\pi\al)u^\al+1}, \quad u>0,
\end{equation}
is a p.d.f. These expressions were not used in our approach.
\item Using \eqref{Form:LaplaceDal}, we see that Theorem \ref{Th:MittagLeffler2GGC} could be restated as follows:  for $\alpha \in (0,1)\cup (1,2]$, we have
\begin{equation}\label{gg}
\lambda \mapsto \frac{\al}{|\al-1|} \frac{\la^{\al-1}-1}{\la^\al-1}\; \II_{\la\ne1}+  \;\II_{\la=1}
    \in \hcm \Longleftrightarrow \al\geq \frac{1}{2}.
\end{equation}
In \cite[Theorem~5.7.1 ]{B}, we have the following computation, valid  for all $c\in (0,1)$:
$$\er\left[\Exp{-\la \frac{\ga_1 \times\ga_{1-c}}{\ga_c}}\right] = \frac{1-\la^c}{1-\la}\text{ if }\la\ne1\quad \et \quad\er\left[\Exp{- \frac{\ga_1 \times\ga_{1-c}}{\ga_c}}\right] = c.$$
Note that the latter function is $\hcm$, since it is the Laplace transform of the product and quotient of independent $\HCM$ random variables. Using property \eqref{propel}, we see that
$$\lambda \mapsto \frac{1}{c}   \frac{1-\la^{\theta c}}{1-\la^\theta} \;\II_{\la\ne1} +\II_{\la=1} \in \hcm, \quad \forall \theta\in [-1,1].$$
Taking  $c=|\al-1|/\al$ for $\alpha \in [1/2,1)\cup (1,2]$, we retrieve the ``if part'' in \eqref{gg}, and also point 1) in Theorem~\ref{Th:MittagLeffler2GGC}.  When $\al>2$, the function in \eqref{gg} does not extend  anymore to an analytic function on $\Cset\setminus(-\infty,0]$, and by \cite[(ix) p. ~68]{B}, it follows that this function cannot be $\hcm$.
\end{enumerate}
\subsection{Comments on Theorem \ref{sin}.}
In \eqref{wella}, the l.h.s. term is certainly the Laplace transform of a distribution on the positive line, whereas the r.h.s. term is the imaginary part of the Laplace transform (with complex argument) of a signed function. We were not able to invert \eqref{wella} to get the explicit distribution of $\dac^{\;[t]}$ with elementary computations. On the other hand, formula  \cite[FI II 812, BI (361)(9), pp 498]{grad} asserts that
$$\int_{0}^{\infty}  x^{\mu-1} \;e^{-a x}\;\sin(bx) dx=\frac{\Gamma(\mu)}{(a^2+b^2)^{\frac{\mu}{2}}} \sin\left(\mu \;\arctan \frac{a}{b}\right), \quad \Rel(\mu)>-1, \; \Rel(a)>|\Ima(b)|.$$
Using \eqref{wella}  and last formula, we may write: for all  $\lambda\in (-t,t)$,
\begin{eqnarray}
\er[\dac^{t+\lambda}]&=&\frac{1}{\Gamma(t -\lambda)} \int_{0}^{\infty} u^{t -\lambda-1}\; \er[\dac^{2t}\;e^{-u \dac}]\;du\nonumber\\
&=& \frac{t\Gamma(2t)}{\sin(\pi\al)\;\Gamma(t -\lambda)}\; \int_{0}^{\infty} u^{-\lambda-1}\; \er\left[\frac{e^{-\cos(\pi\al)(u/\mathbb{S}_t)^t}}{\mathbb{S}_t^t} \; \sin\big(\sin(\pi\al)\; (u /\mathbb{S}_t)^t\big)\right]\; du\nonumber\\
&=& \frac{\Gamma(2t)}{\sin(\pi\al)\;\Gamma(t -\lambda)}\; \int_{0}^{\infty}
\er\left[\mathbb{S}_t^{-t-\lambda}\right] \;v^{-\frac{\lambda}{t}-1}\; e^{-\cos(\pi\al)(v)}\; \sin\big(\sin(\pi\al)\; v\big)\; dv\nonumber\\
&=&\frac{\Gamma(2t)}{\sin(\pi\al)\;\Gamma(t-\lambda)}\;\frac{\Gamma\left(2+\frac{\lambda}{t}\right)}{\Gamma(1+t+\lambda)}\;
 \Gamma\left(-\frac{\lambda}{t}\right) \sin\left(-\frac{\lambda}{t} \;\pi\al\right)\nonumber\\
&=& \frac{\Gamma(2t)}{\sin(\pi\al)\;\Gamma(t-\lambda)}\;\frac{\Gamma\left(1+\frac{\lambda}{t}\right)}{\Gamma(t+\lambda)}\;
 \Gamma\left(1-\frac{\lambda}{t}\right) \frac{\sin\left(\frac{\lambda}{t} \;\pi\al\right)}{\lambda}.\label{div}
\end{eqnarray}
Letting $\lambda \to 0$, we get
$$\er[\dac^{t}]= \frac{\pi\;\al\;\Gamma(2t)}{t\;\sin(\pi \al)\; \Gamma(t)^2}.$$
Then, dividing  \eqref{div} by the latter, we retrieve another computation of the  complicate gamma ratio  form of the Mellin transform \eqref{carcd}.  After searching in several books specialized in integral representations, such as \cite{grad}, we were unsuccessful in finding this complicated expression. Finally, using identity \eqref{ideng}, the fact that $\dac^{\;[t]}\sim 1/\dac^{\;[t]}$, and the beta prime p.d.f. of $\ga_t/\ga_t'$ explicitly given by $\Gamma(2t)\;\Gamma(t)^{-2}\; x^{t-1}\;(1+x)^{-2t}, \; x>0$, we get the alternative representation
\begin{equation*}\label{ftal0}
f_{\tal^{1/t}}(x)=\frac{t \sin(\pi \alpha)}{\pi \alpha}\frac{x^{t-1}}{1+2 \;\cos(\pi \alpha) \;x^t + x^{2t}}= \frac{\Gamma(2t)}{\Gamma(t)^{2}}\; x^{t-1}\;\er\left[\frac{\left(\dac^{\;[t]}\right)^t}{\left(x+\dac^{\;[t]}\right)^{2t}}\right], \;\;\; x>0,
\end{equation*}
where that parameters are as in Theorem \ref{sin}. In particular,
\begin{equation}\label{party}
x\mapsto \frac{1}{1+2 \;\cos(\pi \alpha) \;x^t + x^{2t}}\in \ST_{2t}.
\end{equation}
\subsection{Comments on Corollary \ref{corg}.}
Point 2) in Corollary \eqref{corg} completes \eqref{bobo} and comforts conjecture \eqref{conju}. We aim to deepen the investigation of this conjecture in a forthcoming work.
\subsection{Some account of stable distributions.}
We will need the several identities for the  positive stable distributions.  Shanbhag and  Sreehari  \cite[Theorem 1]{K} exhibited following gamma-mixture identity:
\begin{equation}\label{ss}
\ga_r^{1/\alpha}\simdis \frac{\ga_{\alpha r}}{\sta_{\alpha}^{\,[-\alpha r]}},\quad r>0, \;\; 1>\alpha>0.
\end{equation}
The following one is more classical:
\begin{equation} \label{facto}
\ga_1 \simdis  \left(\frac{\ga_1}{\sal} \right)^\alpha, \quad\mbox{where in the l.h.s, the exponentially distributed r.v. $\ga_1$  is independent of $\sal$,}
\end{equation}
and can be easily seen from
$$\pr (\exs^{1/\alpha} > \lambda) =\pr (\exs > \lambda^\alpha) =e^{-\lambda^\alpha} = \er[e^{-\lambda \mathbb{S}_\alpha}]= \pr (\exs > \lambda \sal)= \pr \left(\frac{\exs}{\sal}  > \lambda\right),\quad \lambda\geq 0.$$
Note that the last identity yields the Mellin transform computation in \eqref{moms}, and that the explicit expression of the p.d.f. of $\tal$ is easily retrieved through the expression of the Mellin transform of $\tal$ obtained from \eqref{moms}, the Euler's reflection formula for the gamma function, and the residue theorem. Additionally, if $|s|<1$, then
\begin{equation}\label{momt}
\er\left[ \tal^{s} \right]=\frac{\Ga(1-s)\;\Ga(1+s)}{\Ga(1-\alpha s)\;\Ga(1+\alpha s)}=\frac{\sin(\pi \alpha s )}{\sin(\pi s)}= \frac{\sin(\pi \alpha )}{\pi}
\int_0^\infty   \frac{x^s}{1+ 2 \cos(\pi \alpha)x +x^2}\;dx,
\end{equation}
and the injectivity of the Mellin transform provides the p.d.f. of $\tal$ in \eqref{stst}. Also note that a possible approach to conjecture \eqref{conju} is Kanter's factorization that can be found in \cite[Corollary 4.1]{K}:
\begin{equation}\label{kan}
\frac{1}{\sal^{\al/(1-\al)}}\simdis  \frac{\ga_1}{s_{\al}^{1/(1-\alpha)}(\mathbb{U})},
\end{equation}
where $\mathbb{U}$ denotes an uniform random variable on $(0,1)$ independent of $\ga_1$ and $s_{\al}$ is the function defined on $(0,1)$ by
$$s_{\al}(u)=\frac{\sin^\al\big(\pi\al u\big)\,\sin^{1-\al}\big(\pi(1-\al)u\big)}{\sin(\pi u)}, \quad u\in (0,1).$$
Combining identities \eqref{ss} and \eqref{kan} with the beta-gamma algebra \eqref{alge}, we can complete \eqref{boto} by the following: for $\al\in(0,1)$, it holds that
\begin{equation}\label{kon}
0<t\leq   1-\al \Longleftrightarrow \sal^{-\frac{\al}{t}} \sim\ID  \Longleftrightarrow \sal^{-\frac{\al}{t}} \mbox{ is a $\ga_1$-mixture} \Longrightarrow \tal^{\frac{1}{t}}\mbox{ is a $\ga_1$-mixture} \Longrightarrow
 \tal^{\frac{1}{t}}\sim \ID.
\end{equation}
Thus, \eqref{boto}, \eqref{bobo} and \eqref{kon} appear to give a certain credit to conjecture \eqref{conju}.

\subsection{Some account of Thorin and complete Bernstein functions.}
Note that function $\phi$ belongs to $\TB$ (respectively $\CB$) if $\phi\in \BF$ if its L\'evy measure $\Pi$  in \eqref{BF} has a density function of the form $\Lap_U(x)/x$ (respectively $\Lap_V(x)$), $x>0$,  where the  $U$ and $V$ and are  positive  measures on $(0,\infty)$  satisfying \eqref{lolo}. Both classes $\CB$ and $\TB$ are convex cones that are closed under pointwise limits; the class $\CB$ is stable by composition, and by  \cite[Theorem 8.4]{SSV}, for a Thorin Bernstein function $\phi$, we have
\begin{equation}\label{compotb}
\varphi \circ \phi \in \TB, \quad \forall \varphi \in \TB \Longleftrightarrow \frac{\phi'}{\phi} \in \ST_1.
\end{equation}
It is immediate that
\begin{equation}\label{ps}
\phi \in \CB  \Longleftrightarrow \lambda \mapsto \phi\left(\frac{1}{\lambda }\right) \;\in
\ST_1  \Longleftrightarrow \quad \lambda \mapsto\frac{\phi(\lambda)}{\lambda } \;\in
\ST_1 \,.
\end{equation}
The following fact is much less evident. By \cite[Theorem 7.3]{SSV} , we also have
\begin{equation}\label{toto1}
\phi \;\in \CB \;   \Longleftrightarrow \; \frac{1}{\phi }\;\in \ST_1 \,.
\end{equation}
An  important representation for the logarithmic derivative of $\phi \in \CB$ is provided by \cite[Theorem 6.17]{SSV}:  there exists  (a unique) pair $\gamma \in \Rset$ and a measurable function $\eta:\oi\to [0,1]$ such that
\begin{equation}\label{repstielcbf'f}
\frac{\phi'(\lambda)}{\phi(\lambda)}=\int_{\oi}\frac{\eta(u)}{(\lambda +u)^2} du, \quad \lambda >0.
\end{equation}
For instance, the Thorin Bernstein function $\phi_\alpha(\lambda)=\lambda^\alpha, \;\alpha \in (0,1)$, has   $\eta$-function equal to $\alpha \II_{u>0}$.
\section{The proofs.}\label{theproofs}
\begin{proofof}{\bf Theorem \ref{mi1}.}  1) From \cite[(3.4) and (3.6)]{SIM3}, we have
$$e^\la-\al  E_\al(\la^\al)=\frac{\al\; \sin(\pi\al)}{\pi}\int_0^\infty e^{-\la t}\frac{t^{\al-1}}{t^{2\al}-2\cos(\pi\al)t^\al+1}dx, \; \;\la\ge0.$$
Aftre some identification in \eqref{ftal},  we get the Laplace representation  \eqref{Mittag1}:
$$e^\la-\al  E_\al(\la^\al) =(1-\alpha)  \begin{cases}
\er\left[e^{-\la \T_{1-\al}^{1/\al}}\right], \; & \text{if $\al\in(0,1)$,}\\
&\\
\er\left[e^{-\la \T_{\al-1}^{1/\al}} \right], \; & \text{if $\al\in(1,2]$}.
\end{cases}
$$
2) a) We have
\begin{equation}\label{notia}
1- \alpha \;e^{-\lambda} E_\alpha (\lambda^\alpha) =(1-\alpha)  \; \er[e^{-\lambda {\mathbb{R}_\alpha}}], \quad e^{-\l} E_\alpha (\lambda^\alpha)-1= \frac{1-\alpha}{\alpha}  \;\big(1-\er[e^{-\lambda {\mathbb{R}_\alpha}}]\big), \quad \lambda \geq 0.
\end{equation}
Let $\mathbb{R}_n^{(\al)}$ be the random walk generated by $\mathbb{R}_\alpha$  given by  \eqref{tal}, i.e. $\mathbb{R}^{(\al)}_0=0,\;\;\;\mathbb{R}_n^{(\al)}=\mathbb{R}_{1,\al}+\ldots+\mathbb{R}_{n,\al},\; n\geq 1,$ where $\mathbb{R}_{1,\al},\ldots, \mathbb{R}_{n,\al}$, are $n$ independent copies of $\;\mathbb{R}_\al$. Let $\mathbb{Z}_\al$ be geometrically distributed on $\Nset$ with parameter $\alpha$, and independent of the sequence ${(\mathbb{R}_n^{(\al)})}_{n\geq 0}$.  Using  \eqref{tal}, \eqref{Mittag1} and the subordinated r.v. $\mathbb{R}^{(\al)}_{\mathbb{Z}_\al}$, we have the Laplace representation
\begin{eqnarray*}
\frac{e^{\lambda}}{ E_\alpha(\lambda^\alpha)} &=& \frac{\alpha}{ 1- (1-\alpha)\er[e^{-\lambda \mathbb{R}_\al}] }= \sum_{n=0}^\infty \alpha\;(1-\alpha)^n \er[e^{-\lambda \mathbb{R}_\al}]^n
=   \sum_{n=0}^\infty \pr(\mathbb{Z}_\al=n) \;\er[e^{-\lambda \mathbb{R}^{(\al)}_n}]
= \er[e^{-\lambda  \mathbb{R}^{(\al)}_{\mathbb{Z}_\al} }], \quad \lambda \geq 0.
\end{eqnarray*}
2) b) From Kanter's factorization \eqref{kan}, we see  that $\mathbb{S}_{1-\al}^{-(1-\al)/\al}$ has a completely monotone density, and so does  $\T_{1-\al}^{1/\al}$. We deduce that  $F_\al$ is a Stieltjes function.  The equivalence between  $F_\alpha \in \hcm$, and  $1/2\leq \alpha <1$, is due to \eqref{bobo}.  The function $-F_\al/F_\al$ is a Stieltjes because $-\log F_\al \in \TB$, and because of \eqref{eqstielcbt}. The function $1- \alpha \;e^{-\lambda} E_\alpha (\lambda^\alpha)$  is $\hcm$ due to \eqref{noti} and \eqref{notia}.\\
3) a) We do as in 2) a). \\
3) b) By \eqref{bobo}, we know that if $t>2(2-\alpha),$ then the distribution of  $\T_{\al-1}^{1/t}$ is not $\GGC$. Since $\alpha>2(2-\alpha)$ in case $\alpha >4/3$, then $F_\alpha$ could not be $\hcm$.\\
3) c) For the last assertion,  use \eqref{xx}, \eqref{stal121} and \eqref{bobo}, and  observe that $0<\gamma\leq  (2-\alpha)/\alpha<1$. Then
\begin{eqnarray*}
1<\alpha \leq \frac{3}{2} &\Longleftrightarrow &   0<\alpha-1<\frac{1}{2}\Longrightarrow \mathbb{S}_\gamma \;\mathbb{T}_{\alpha -1}^{1/\gamma\alpha} \sim \GGC\\
   &\Longrightarrow&  F_\alpha(\lambda^\gamma)=\er\left[e^{-\lambda^\gamma {\mathbb{T}_{\alpha -1}}^{1/\alpha}}\right]= \er\left[e^{-\lambda \mathbb{S}_\gamma \;\mathbb{T}_{\alpha -1}^{1/\gamma\alpha}}\right]\in \hcm.
\end{eqnarray*}
As in 2), the mixture property \eqref{kon}, entails $\mathbb{T}_{\alpha -1}^{1/\gamma\alpha}$ has a completely monotone density, hence so does $\mathbb{S}_\gamma \;\mathbb{T}_{\alpha -1}^{1/\gamma\alpha}$.
The latter proves that $F_\al(\la^\gamma)$ is Stieltjes; for it logarithmic derivative, conclude as in the end of the proof of point 2) b).
\end{proofof}
\medskip

\begin{proofof}{\bf Theorem \ref{Th:MittagLeffler2GGC}.}
\noindent 1) Assume $\al\in(1,2]$.  The function $\Lap_{D_2}(\la)= (\la+1)^{-1}$ is clearly an $\hcm$ function, and by \eqref{gcgc}, we deduce $D_2 \sim\GGC$.  If $\al\in(1,2)$, then the function
\begin{equation}\label{abo}
\frac{1+x}{x^{2\al}-2\cos(\pi\al)x^\al +1},\quad x>0,
\end{equation}
is locally increasing in a neighborhood of $0+$ and cannot be completely monotone. Suppose that $G_\al \in \hcm$. By   \cite[Theorem 5.4.1]{B}, the function $G_\al$ is then the Laplace transform of a $\GGC$ and by   \cite[Theorem 4.1.1 and Theorem 4.1.4]{B}, the Thorin mass in \eqref{abov}  equals $\al$, thus the function in \eqref{abo} would be completely monotone, a contradiction. We deduce that $G_\al$ is not $\hcm$ when $\al\in(1,2]$. To show that $G_\al \sim \GGC$,  it suffices to show that $G_{1,\alpha}:=-\log\Lap_{G_\al}$  is a Thorin-Bernstein function (equivalently $G_{1,\alpha}'$ is a Stieltjes function), using formula~(\ref{Form:LaplaceDal}).  The function $G_{1,\alpha}$ (and then $G_{1,\alpha}'$)  extends to  an analytic function on $\Cset\setminus\Rset_-$. We will use the characterization of Stieltjes function given \cite[by Corollary~7.4]{SSV}, namely, we aim to prove that
\begin{equation}\label{im}
\Ima(z)>0\Longrightarrow \Ima(G'(z))<0.
\end{equation}
Observe that $\Ima(G_{1,\alpha}')$ is harmonic on the upper half-plane as the imaginary part of an analytic function. Moreover, $G_{1,\alpha}'(z)\to 0$ uniformly as $|z|\to+\infty$. Then, from a compactness argument and from the minimum principle, it suffices to show that  $$\limsup_{z\to0,~\Ima(z)>0}\Ima(G_{1,\alpha}'(z))\le0.$$ Elementary computations give,
$$G_{1,\alpha}'(z)=\frac{\al z^{\al-1}}{z^\al-1}-\frac{(\al-1)z^{\al-2}}{z^{\al-1}-1}, \quad z\in\Cset\setminus\Rset_-,$$
and for all $x>0$, we have
\begin{align*}
G_{1,\alpha}'(-x^+)&:=\lim_{\scriptsize \begin{matrix}
z\to -x\\\Ima(z)>0\end{matrix}}
G_{1,\alpha}'(z)  = e^{i\pi\al}x^{\al-2}\pa{\frac{\al-1}{x^{\al-1}e^{i\pi\al}+1}-\frac{\al x}{x^\al e^{i\pi\al}-1}}\\
& =- x^{\al-2}\,\frac{x^\al e^{i2\pi\al}+(\al-1)e^{i\pi\al}+\al x e^{i\pi\al}}{x^{2\al-1}e^{i2\pi\al}+(x^\al-x^{\al-1})e^{i\pi\al}-1},\\
\Ima(G_{1,\alpha}'(-x^+))&=\sin(\pi\al)\;x^{\al-2}\;\frac{A_{1,\al}(x)}{B_{1,\al}(x)},
\end{align*}
where
$$A_{1,\al}(x)=(\al-1)x^{2\al}+\al x^{2\al-1}+2\cos(\pi\al)x^\al+\al x+\al-1,  \quad
B_{1,\al}(x)=\left| x^{2\al-1} e^{i2\pi\al}+(x^\al-x^{\al-1})e^{i\pi\al}-1 \right|^2.$$
The function $A_{1,\al}$ is a positive  on $[0,\infty)$. Indeed,
\begin{align*}
A_{1,\al}'(x) & = \al\pab{ 2(\al-1)x^{2\al-1}+(2\al-1)x^{2\al-2}+2\cos(\pi\al)x^{\al-1}+1}\\
& \ge \al\pa{x^{2\al-2}-2x^{\al-1}+1}=\al\pa{x^{\al-1}-1}^2\ge0,
\end{align*}
thus, $A_{1,\al}$ is a non-decreasing function on $[0,\infty)$, and since $A_{1,\al}(0)=\al-1>0$, we deduce that $A_{1,\al}$ is positive and then $\Ima(G_{1,\alpha}'(-x^+))\le0$ for all $x>0$. Next, since $-1<\al-2<0$, then
$$\limsup_{\scriptsize \begin{matrix}
z\to 0\\\Ima(z)>0\end{matrix}}\Ima(G_{1,\alpha}'(z))=\limsup_{\scriptsize \begin{matrix}
z\to 0\\\Ima(z)>0\end{matrix}}\Im\pa{(\al-1)z^{\al-2}}=(\al-1)\limsup_{\scriptsize \begin{matrix}
r\to 0,~r>0\\\theta\in(0,1)\end{matrix}} r^{\al-2}\sin(\pi(\al-2)\theta)=0,$$
Finally, $\Ima(G_{1,\alpha}'(x))=0$, for all$\, x>0$.
All in all, we have proved that the function $G_{1,\alpha}'$ is a Stieltjes, hence  $G_{1,\alpha}$ is  Thorin-Bernstein.\\
2) Assume $\al\in(0,1)$. Using \cite[Theorem 5.4.1]{B} again and the definition of a widened $\GGC$ in \cite[Section 3.5]{B}, we obtain
$$-G_\al\;\; \mbox{is a widened} \;\GGC\Longleftrightarrow \la\mapsto \Lap(-G_\al)(\lambda)=\frac{1-\la^{\al-1}}{\la^\al-1}\in \hcm\Longleftrightarrow\la\mapsto \frac{\la^{1-\al}-1}{\la^\al-1}\in \hcm,$$
and it suffices to prove that the function
$$G_{2,\alpha}(\la):=-\log\pa{\frac{\la^{1-\al}-1}{\la^\al-1}}, \quad \lambda \geq 0,$$
is  Thorin-Bernstein (equivalently $G_{2,\alpha}'$ is a Stieltjes function) if and only if $\al\ge1/2$.
Elementary computations give,
$$G_{2,\alpha}'(z)=\frac{\al z^{\al-1}}{z^\al-1}-\frac{(1-\al)z^{-\al}}{z^{1-\al}-1},\quad  z\in\Cset\setminus(-\infty,0],$$
and then,
$$
\Ima(G_{2,\alpha}'(-x^+))=-\sin(\pi\al)\frac{A_{2,\al}(x)}{B_{2,\al}(x)}, \quad  x>0,$$
where
$$A_{2,\al}(x)=(1-\al)(x^\al+x^{-\al})-\al (x^{1-\al}+ x^{\al-1})-2\cos(\pi\al)
\quad \mbox{and}\quad  B_{2,\al}(x)=\left|x-1+x^\al e^{i\pi\al}-x^{1-\al}e^{-i\pi\al}\right|^2.$$
We proceed as we did for $G_{1,\alpha}'$ to show that $G_{2,\alpha}'$ is a Stieltjes function.  For this, we have to check the sign of $A_{2,\al}$.  The function $H_{1/2}=0$ is a trivial Thorin-Bernstein function. Since for $\al\ne1/2$
$$A_{2,\al}(1)=2(1-2\al-\cos(\pi\al))~ \begin{cases}
<0,& \text{ if }\al \in(0,1/2)\\
>0, & \text{ if }\al \in(1/2,1),
\end{cases}$$
we see that  $\al\in(1/2,1)$ is a necessary condition for $A_{2,\al}$ to be positive. Next, let
$$C_\al(u)=A_{2,\al}(e^u)=2\pab{(1-\al)\cosh(\al u)-\al\cosh(\al u)-\cos(\pi\al)}, \quad u\in\Rset.$$
With the expression
$$C_\al'(u)=2\al(1-\al)\pab{\sinh(\al u)-\sinh((1-\al)u)},$$
one deduces that  if $\al\in(1/2,1)$, then the function $A_{2,\al}$  decreases on $(0,1)$ and increases on $(1,+\infty)$, thus  $A_{2,\al}(x)\ge A_{2,\al}(1)>0$ for all $\;x>0$. Finally, we deduce that if $\al\in(1/2,1)$, then
$$\limsup_{\scriptsize \begin{matrix}
z\to 0\\\Ima(z)>0\end{matrix}}\Ima(G_{2,\alpha}'(z))= \limsup_{\scriptsize \begin{matrix}
r\to 0,~r>0\\\theta\in(0,1)\end{matrix}} \sin(\pi\al\theta)(\al r^{\al-1}-(1-\al)e^{-\al}) =0,$$
and this shows  that $G_{2,\alpha}'$ is a Stieltjes function.\\
The $\GGC$  property of $\ga_1/\mathbb{R}_\al$ is straightforward thanks to \eqref{hh}.
\end{proofof}
\medskip

\begin{proofof}{\bf Corollary \ref{to1}.} If $\alpha=2$, then $E_2(\lambda)=\cosh(\sqrt{\lambda})$, ${(X_t^{(2)})}_{t\geq 0}$ is Brownian motion, and $\tau_1\simdis \mathbb{S}_{1/2}$. The rest of the assertions is straightforward due to the convention  $\mathbb{T}_{0}=1$. We then treat the case  $\alpha\in (1,2)$.\\
1) Using \cite[(1.4)]{SIM2} and the definition of $F_\alpha$ in \eqref{Mittag1}, one has the Laplace representation
\begin{equation}\label{dou}
\er[e^{-\lambda^\alpha\; \tau_1}]=(1-\alpha)F_\alpha'(\lambda)=H_\alpha(\lambda^\alpha),
\end{equation}
and this gives the Laplace transform representation \eqref{galga}. \\
2) From \eqref{Mittag1} and \eqref{dou}, we obtain
\begin{equation}\label{obo}
\er[e^{-\lambda\; \tau_1}] =(\alpha-1)\;\er\left[{\talm}^{1/\alpha}  \;e^{-\lambda^{1/\alpha} {\talm}^{1/\alpha} }\right]= (\alpha-1)\;\er\left[{\talm}^{1/\alpha}  \;e^{-\lambda \;\mathbb{S}_{1/\alpha} \talm}\right],
\end{equation}
where $\mathbb{S}_{1/\alpha}$ is a positive r.v. independent of $\tal$. Further, \eqref{obo} or  \eqref{d0} imply
\begin{equation*}\label{alta}
1=(\alpha-1)\;\er[{\talm}^{1/\alpha}],
\end{equation*}
hence
$$\er[e^{-\lambda\; \tau_1}] = \frac{\er[{\talm}^{1/\alpha}  \;e^{-\lambda \;\mathbb{S}_{1/\alpha} \talm}]}{\er[{\talm}^{1/\alpha}]}=\er\left[e^{-\lambda \;\mathbb{S}_{1/\alpha} (\talm)^{\;[1/\alpha]}}\right],$$
which shows the first factorization in \eqref{tau2}. The second factorization in \eqref{tau2} is due to \eqref{facto} and to the second identity in \eqref{biaiso} applied to $\talm$. \\
3) The second factorization in \eqref{tau2} gives
$$\left(\frac{\ga_1}{\tau_1}\right)^{1/(2-\alpha)} \simdis  \ga_1^{\alpha/(2-\alpha)}\; \left(\talm^{1/(2-\alpha)}\right)^{\;[(2-\alpha)/\alpha]}.$$
Since $q=\alpha/(2-\alpha)>1$, then the power property in \eqref{xx} ensures that $\ga_1^{\alpha/(2-\alpha)} \sim \HCM$. Property \eqref{bobo} with $t=(2-\alpha)$ asserts that $\talm^{1/(2-\alpha)}\sim \HCM$, and by property \eqref{prov} we deduce that $\left(\talm^{1/(2-\alpha)}\right)^{\;[(2-\alpha)/\alpha]}\sim \HCM.$
Finally, the stability property by independent products in \eqref{xx} yields  $(\ga_1/\tau_1)^{1/(2-\alpha)}\sim \HCM$.\\
4) We know that $\mathbb{S}_{1/\alpha}\sim \GGC$. Then
$$0<\alpha-1 \leq  \frac{1}{2} \stackrel{\eqref{bobo}}{\Longrightarrow}  \talm\sim\HCM \stackrel{\eqref{prov}}{\Longrightarrow}(\talm)^{\;[1/\alpha]}\sim\HCM\stackrel{\eqref{xx}}{\Longrightarrow}
\tau_1 \simdis \mathbb{S}_{1/\alpha}\,(\talm)^{\;[1/\alpha]} \sim\GGC,$$
the $\HCM$ property of $H_\alpha$ follows from \eqref{gcgc}.
\end{proofof}

\medskip
\begin{proofof}{\bf Proposition \ref{hcmlog}.} 1) By \eqref{compotb}, \eqref{ps} and \eqref{toto1},  we have
\begin{eqnarray*}\phi_{c,t} \in \TB &\Longleftrightarrow& \lambda \mapsto \phi_{c,t}'(\lambda)= \frac{2t}{\lambda} \, \frac{c \lambda^{t} +\lambda^{2t}}{1+2 c \lambda^t + \lambda^{2t}} \in \ST_1 \Longleftrightarrow   \lambda \mapsto \frac{c \lambda^{t} +  \lambda^{2t}}{1+2 c \lambda^t + \lambda^{2t}} \in \CB \\
&\Longleftrightarrow& \lambda \mapsto \frac{1+2 c \lambda^t + \lambda^{2t}}{c \lambda^{t} +  \lambda^{2t}}=1+\frac{1+c\lambda^t }{c\lambda^{t} + \lambda^{2t}} \in \ST_1  \Longleftrightarrow  \lambda \mapsto \frac{1+ c \lambda^t }{c \lambda^{t} +  \lambda^{2t}} \in \ST_1 \\
& \Longleftrightarrow& \lambda \mapsto  \varphi_{c,t}(\lambda):=\frac{c \lambda^{t}+\lambda^{2t}}{1+ c \lambda^t }\in\CB.
\end{eqnarray*}
1- \underline{The case  $t = 1$}. We have $\varphi_{1,1}(\lambda)=1$ and the function $\varphi_{c,1}(\lambda)/\lambda=(c+\lambda)/(1+c\lambda)$ is not completely monotone if $c\neq 1$, because its derivative $c(1-c)(1-\lambda)/(1+c\lambda)^2$ has a change sign. Thus, $\varphi_{c,1} \notin \BF$ if $c\neq 1$.\\

\noindent  2- \underline{The case  $c = 1$}. Trivially, $\varphi_{1,t}(\lambda)= \lambda^{t}\in \CB\Longleftrightarrow t\leq 1$.\\

\noindent 3- \underline{The case  $c \neq1$}. The function $\varphi_{c,t}$ is $\CB$ if, and only if, $\varphi_{c,t}$ satisfies \eqref{repstielcbf'f}. We then study the logarithmic derivative of $\varphi_{c,t}$:
\begin{eqnarray*}\frac{\varphi_{c,t}'}{\varphi_{c,t}}(\lambda)  &=& \frac{t}{\lambda} \,\left[\frac{c \lambda^{t} +  2 \lambda^{2t}}{c \lambda^{t} +  \lambda^{2t}} - \frac{c \lambda^{t} }{1+ c \lambda^{t}} \right]=\frac{t}{\lambda} \,\left[\frac{ \lambda^{2t}}{c \lambda^{t} + \lambda^{2t}} +1 - \frac{c \lambda^{t} }{1+ c \lambda^{t}} \right]=\frac{t}{\lambda} \,\left[\frac{\lambda^t}{c +\lambda^{t}} + \frac{1 }{1+ c \lambda^{t}} \right].
\end{eqnarray*}

3-a) \underline{The case  $c \neq1,\, t>1$}.  Here, we have a first conclusion:
if $t>1$, then $\varphi_{c,t}$ does not belong to $\CB$, and not even to $\BF$,  since
$$    \frac{\varphi_{c,t}(\lambda)}{\lambda}\sim  \lambda^{t-1} \to  +\infty, \quad \mbox{when} \;\lambda \to +\infty\,.$$

3-b) \underline{The case  $c \neq1,\, 0<t<1$}. Recall the r.v. $\Tt$ given in (\ref{stst}). Observe that if  $\ga_1$ is standard exponentially  distributed and is independent of $\Tt$, then
$$\frac{1}{1+\lambda^t}= \er[e^{-\lambda \ga_{1}\Tt^{1/t}}], \quad \frac{\lambda^t}{1+\lambda^t}= 1-\er[e^{-\lambda \ga_{1}\Tt^{1/t}}].$$
Using the fact that used that $\Tt^{1/t}\simdis \Tt^{-1/t}$ and the latter, we obtain
\begin{eqnarray*}\frac{1}{\lambda}  \frac{1 }{1+ c \lambda^{t}} &=& \frac{1}{\lambda}  \frac{1 }{1+ (c^{1/t} \lambda)^{t}}= \frac{1}{\lambda} \er \left[e^{-\lambda \ga_{1}(c \Tt)^{1/t}}\right]= \int_0^\infty e^{-\lambda x}  \,  \pr(\ga_{1}(c \Tt)^{1/t}  \leq x)\, dx\\
&=& \int_0^\infty e^{-\lambda x} \,x  \, \frac{1- \er[e^{ - x(\Tt/c)^{1/t} }]}{x}  dx=\int_0^\infty e^{-\lambda x} \,x \, \Lap_f(x)\,dx=\int_0^\infty \frac{f(u)}{(\lambda +u)^2} du,
\end{eqnarray*}
where,
$$f(u)= \pr(\Tt^{1/t} > u c^{1/t})= \pr(\Tt  > c u^t), \quad u \geq 0\,.$$
Similarly, we have
\begin{eqnarray*}\frac{1}{\lambda}  \frac{\lambda^{t} }{\lambda^{t}+ c } &=& \frac{1}{\lambda}  \frac{(\lambda/c^{1/t})^t }{1+ (\lambda/c^{1/t})^{t}}= \frac{1}{\lambda} \left (1- \er \left[e^{-\lambda \ga_{1} (\Tt/c)^{1/t}}\right]\right) =\int_0^\infty e^{-\lambda x}  \,\pr(\ga_{1}(\Tt/ c)^{1/t}  > x)\,dx\\
&=&\int_0^\infty e^{-\lambda x} \,x\,\frac{\er\left[ e^{- x (c\Tt)^{1/t}} \right]}{x}dx=\int_0^\infty e^{-\lambda x} \,x \, \Lap_g(x)  \, dx =\int_0^\infty \frac{g(u)}{(\lambda +u)^2} du,
\end{eqnarray*}
where
$$g(u)= \pr\big (c \Tt)^{1/t} < u \big) =\pr(\Tt  < u^t/c), \quad u \geq 0\,.$$
We finally get the representation of logarithmic derivative of $\varphi_{c,t}$:
$$\frac{\varphi_{c,t}'}{\varphi_{c,t}}(\lambda)= \int_0^\infty \frac{\eta_t(u)}{(\lambda +u)^2} du, \qquad \eta_t(u):= t\,\big[ f(u)+g(u)\big] = t\,\left[ \pr(\Tt  > c u^t)+\pr(\Tt  < u^t/c)\right].$$
Due to (\ref{repstielcbf'f}), we can now assert that
\begin{equation}\label{pe}
\phi_{c,t} \in \TB \Longleftrightarrow \eta_{c,t}(u) \in [0,1], \quad \forall u>0\,.
\end{equation}
We arrive to the last step. Since $\Tt$ has the explicit density
$${f}_{\scriptsize \Tt} (u) = \frac{\sin(\pi t)}{\pi t \big(1+2 \cos(\pi t) u+ u^2\big)},\quad u>0\,,$$
then the derivative of $\widetilde{\eta}_{c,t}(u):=\eta_{c,t}(u^t)$ is
\begin{eqnarray*}
\widetilde{\eta}_{c,t}'(u)&=& t\,\left[\frac{1}{c} \,{f}_{\scriptsize \Tt}\left(\frac{u}{c}\right)- c {f}_{\scriptsize \Tt}(c u)\right]=\frac{\sin(\pi t) \,c}{\pi}\,\left[ \frac{1}{c^2 +2 c \cos(\pi t) u +u^2}-  \frac{1}{1+2 c \cos(\pi t) u +c^2\,u^2} \right] ,
\end{eqnarray*}
and finally,
\begin{equation}\label{etap}
\widetilde{\eta}_{c,t}'(u) \geq 0 \Longleftrightarrow 1 +2 c \cos(\pi t) u  +c^2\,u^{2} \geq  c^{2} +2 c \cos(\pi t) u +u^2 \Longleftrightarrow 1- c^2 \geq (1- c^2) u^{2} \Longleftrightarrow u\in [1,\infty).
\end{equation}

3-b)(i) \underline{The case  $c > 1, \,0< t< 1$}. By (\ref{etap}), the function $\eta_{c,t}$ decreases on $[0,1]$, then increases on $[1,\infty)$, and its maximum is equal to $t=\eta_{c,t}(0)=\eta_{c,t}(\infty)$. We deduce that (\ref{pe}) is true.

3-b)(ii) \underline{The case  $c < 1, \,0< t<1$}. By (\ref{etap}), the function $\eta_{c,t}$ increases on $[0,1]$, then decreases on $[1,\infty)$, and its maximum is equal to
\begin{eqnarray*}
\eta_{c,t}(1)&=& t\,\left[ \pr(\Tt  > c)+\pr(\Tt  < 1/c)\right] = 2 t\, \pr(\Tt > c)= 2 t\,\int^\infty_c \frac{\sin(\pi t)}{\pi t (1+2 \cos(\pi t) u +u^2)} \\
&=& \frac{2}{\pi} \left[\frac{\pi}{2} - \arctan \frac{c+\cos(\pi t)}{\sin(\pi t)}\right].
\end{eqnarray*}
Note that we have performed the obvious change of variable $u= \sin(\pi t) v -\cos (\pi t)$ to get the last expression. We deduce that (\ref{pe}) is true if, and only if,
 $$\arctan \frac{c+\cos(\pi t)}{\sin(\pi t)} \geq 0 \Longleftrightarrow c+\cos(\pi t) \geq 0\,.$$
2) By 1), we have $(i) \Longleftrightarrow (iv)$. Trivially, $(iv)\Longrightarrow (iii)\Longrightarrow (ii)$. It suffices to check $(ii)\Longrightarrow (1)$. Assume
\begin{equation*}\label{last}
g_{\alpha,t}(\lambda):= \frac{1}{1+2 \cos(\pi \alpha) \lambda^t +  \lambda^{2t}}\in \CM.
\end{equation*}
If $\alpha \in (1/2,1)$, then  the denominator has a sign change, and $g_{\alpha,t}$ could not be $\CM$.  Assume $\alpha \in (0,1/2]$ and $t>1-\alpha$, then  $g_{\alpha,t}$ has two poles  $e^{\pm \ii \pi (1-\alpha)/t}$. In this case $g_{\alpha,t}$ can not extend to an analytic function on $\{z\in\Cset\;/\;\Rel(z)>0\}$, then, it could not be $\CM$.
\end{proofof}
\medskip

\begin{proofof}{\bf Corollary \ref{crt1}.} 1) Using \eqref{mittagt}, note that $$\chi_c(u)=\Rel\left(e^{-\jc\;u}\right)\;\;\;\mbox{and}\;\;\;\mathfrak{C}_{c,t}(x)=\Rel\left(E_t(- \jc\;x^t)\right)=\Rel\left(\er\left[\chi_c\big((x/\sat)^t\big)\right]\right) .$$
2) The equivalences are obtained by the form of the derivative of $\phi_{c,t}$:
\begin{equation}\label{tur}
\phi_{c,t}'(\lambda) = 2\;t\; \frac{c \;\lambda^{t -1}+  \lambda^{2t-1}}{1 + 2 \;c \;\lambda^t+  \lambda^{2t}}= \frac{2t}{\lambda}\;\big(1-s_{c,t}(\lambda)\big),\quad s_{c,t}(\lambda)= \frac{c \;\lambda^t+  1}{1 + 2 \;c \;\lambda^t+  \lambda^{2t}}
\quad \lambda >0.
\end{equation}
By the latter, Proposition \ref{hcmlog}, and by \eqref{eqstielcbt} and \eqref{ps}, we have the equivalences
$$c+\cos(\pi\;t)\geq 0 \Longleftrightarrow\phi_{c,t}\in  \TB\Longleftrightarrow \phi_{c,t}'\in  \ST_1 \Longleftrightarrow 1-s_{c,t}\in  \CB\Longleftrightarrow s_{c,t} \in \ST_1.$$
Using the  representations
$$\chi_c'(u)= -e^{-cu}\left[c\;\cos(\sqrt{1-c^2}\; u)+ \sqrt{1-c^2}\;\sin\big(\sqrt{1-c^2}\; u\big)\right],$$
\begin{equation}\label{xxx}
\frac{x}{1+x^2}=\int_0^\infty e^{-xu} \;\cos(u)du, \;\;\;\; \frac{1}{1+x^2}=\int_0^\infty e^{-xu} \;\sin(u)du , \;\;\;\; x>0,
\end{equation}
and  \eqref{stal}, we may write
\begin{eqnarray*}
s_{c,t}(\lambda)&=& \frac{c \;\lambda^t+  1}{1 + 2 \;c \;\lambda^t+  \lambda^{2t}}= \frac{c}{\sqrt{1-c^2}}  \frac{\frac{\lambda^t+c}{\sqrt{1-c^2}}}{1+\left(\frac{\lambda^t+c}{\sqrt{1-c^2}}\right)^2}+
\frac{1}{1+\left(\frac{\lambda^t+c}{\sqrt{1-c^2}}\right)^2}\\
&=&  \int_0^\infty e^{-(\lambda^t+c)\;u} \;\left[c\;\cos(\sqrt{1-c^2}\; u)+ \sqrt{1-c^2}\;\sin\big(\sqrt{1-c^2}\; u\big)\right]du \\
&=& -\int_0^\infty e^{-\lambda^t \, u} \;\chi_c'(u)du= -\int_0^\infty \er\big[e^{-\lambda \,u^{1/t} \;\sat}\big] \;\chi_c'(u)du.
\end{eqnarray*}
Applying Fubini-Tonelli's theorem, and performing the change of variable $v=u^{1/t} \;\sat$ under the expectation, we obtain the Laplace  transform representation of $s_{c,t}$:
\begin{eqnarray}
s_{c,t}(\lambda)&=& -\int_0^\infty  e^{-\lambda v} \;t\; \;v^{t -1} \er\left[\frac{1}{\sat^{t}} \;\chi_c'\left(\Big(\frac{v}{\sat}\Big)^t\right)\right]\;dv=\int_0^\infty e^{-\lambda v}\left(-\frac{d}{dv} \er\left[\chi_c\left(\Big(\frac{v}{\sat}\Big)^t\right)\right]\right)\;dv \nonumber\\
&=&\int_0^\infty  e^{-\lambda v}  \big(- \mathfrak{C}_{c,t}'(v)\big)\;dv.\label{turr}
\end{eqnarray}
Thus, $s_{c,t}\in \ST_1$, if, and only if, $-\mathfrak{C}_{c,t}'\in \CM$ and $\mathfrak{C}_{c,t}$
(hence decreases from $\mathfrak{C}_{c,t}(0)=1$ to $\mathfrak{C}_{c,t}(\infty)=0$, thus  $\mathfrak{C}_{c,t}\in [0,1]$). The latter is equivalent to $\mathfrak{C}_{c,t}\in \CM$ or to $1-\mathfrak{C}_{c,t}\in \BF$.\\
3) Using \eqref{tur}, \eqref{turr}, and then performing an integration by parts, we obtain the expression
$$\phi_{c,t}'(\lambda)= \frac{2 t}{\lambda} \int_0^\infty  (1-e^{-\lambda v})\big(- \mathfrak{C}_{c,t}'(v)\big)\;dv =2t\int_0^\infty   e^{-\lambda x} \; \mathfrak{C}_{c,t}(x) \;dx.$$
Integrating the last expression from 0 to $\lambda$, we retrieve \eqref{pip}. The  representation of $\mathfrak{C}_{c,t}$ as the Laplace transform or a positive r.v. $\mathbb{E}_{c,t}$ is evident since $\mathfrak{C}_{c,t}$  is completely monotone and $\mathfrak{C}_{c,t}(0)=1$. The last assertion is also evident due to the alternative Frullani integral form of \eqref{pip}.
$$\phi_{c,t}(\lambda)=2 t \;\er
\left[\log\left(1+\frac{\lambda}{\mathbb{E}_{c,t}}\right)\right]=2 t \;\int_{0}^{\infty} \log\left(1+\frac{\lambda}{u}\right) \pr(\mathbb{E}_{c,t}\in du).$$
3) It suffices to write $z=|z|\; \jc$, where $c=\cos(\arg z)$.
\end{proofof}
\medskip

\begin{proofof}{\bf Corollary \ref{hcmcor}.} The first assertion is a straightforward consequence of the Thorin property of the Bernstein function $\phi_{c,t}$ in Proposition \ref{hcmlog}. Indeed, the definition of $\GGC$ distributions gives that
$$\phi_{c,t} \in \TB \Longleftrightarrow x\mapsto e^{-\phi_{c,t}(x)}= \frac{1}{1+2\; c x^t + x^{2t}} \in \hcm.$$
For the second assertion, use the fact that the $\hcm$ class is stable by product, closed by pointwise limits (property (ii) pp. 68 \cite{B}), and property \eqref{propel}.
\end{proofof}
\medskip

\begin{proofof}{\bf Theorem \ref{sin}.} 1) We proceed as in the proof of Corollary \ref{crt1}. Let us define
$$\sigma_c(x):=\Ima \left(\frac{e^{-\jc\, x}-1}{\jc}\right)\quad \mbox{and}\quad \mathfrak{S}_{c,t}(x):= \frac{1}{\sqrt{1-c^2}}\er\left[\sigma_c\left(\Big(\frac{x}{\mathbb{S}_t}\Big)^t\right)\right], \;\;\;x\geq 0.$$
and observe that
$$\sigma_c'(x) = e^{-c x}\sin\big(\sqrt{1-c^2}\, x\big) = \Ima(e^{-\overline{\jc} x}).$$
Using  representations \eqref{lala} and \eqref{xxx},  then performing the change of variable $u=(x/\mathbb{S}_t)^t$ under the expectation, we get
\begin{eqnarray*}
\er[e^{-\lambda \mathbb{X}_{c,t}}]&=&\frac{1}{1+2c\lambda^t+\lambda^{2t}}=
\frac{1}{1-c^2} \frac{1}{1+\left(\frac{\lambda^t+c}{\sqrt{1-c^2}}\right)^2}=
\frac{1}{\sqrt{1-c^2}} \int_0^\infty e^{-(\lambda^t+c)\;u} \; \sin\big(\sqrt{1-c^2}\; u\big) \;du\\
&=& \frac{1}{\sqrt{1-c^2}} \int_0^\infty \er[e^{-\lambda \;u^{1/t} \mathbb{S}_t}]\; \sigma_c'(u) \; du=\frac{t}{\sqrt{1-c^2}} \int_0^\infty e^{-\lambda x}
\;x^{t-1}\;\er\left[ \mathbb{S}_t^{-t} \;\sigma_c'\big((x/\mathbb{S}_t)^t\big)\right]\;dx
\end{eqnarray*}
The latter gives the expression of $f_{_{\mathbb{X}_{c,t}}}(x)$, and by integration over $(0,x]$, we obtain the one of $\mathfrak{S}_{c,t}(x)=\pr(\mathbb{X}_{c,t}\leq x)$ . \\
2) By \eqref{lala}, and Corollary \ref{hcmcor}, recall that the distribution $\mathbb{X}_{c,t}$ is $\GGC$, and is associated with the Thorin Bernstein function  $\phi_{c,t}$ in \eqref{pip} and to the Thorin measure $U_{c,t}$ given by \eqref{repphi3tt}. Using Corollary \ref{crt1}, then \cite[Theorem 4.1.1]{B}, or equivalently \cite[Theorem 4.1.4]{B}, we see that the total mass of $U_{c,t}$ equals
$$U_{c,t}\oi=2 \,t=\sup \,\big\{\,s;\;\lim_{x\to 0+} x^{1-s}f_{\mathbb{X}_{c,t}}(x) =0\,\big\}.$$
and we conclude that  $\mathbb{X}_{c,t}\sim\GGC$  if, and only if, identity \eqref{ideng} holds. The Laplace transform in \eqref{wella} is identified by the equality $\er[\dac^{2t}]=1$ which is obtained by taking $ \lim_{x\to 0+}x^{1-t} f_{\mathbb{X}_{c,t}} (x) $ in \eqref{well}.
By \eqref{fact}, \eqref{well} and \eqref{ideng}, it is straightforward that
\begin{eqnarray*}
\er[\dac^t] \;\er\left[\left(\dac^{\;[t]}\right)^{t}\;e^{-\lambda \dac^{\;[t]}}\right]&=& \er[\dac^{2t}\;e^{-\lambda \dac}]  =
\Gamma(2t)\; \lambda^{1-2t}\; f_{_{\mathbb{X}_{c,t}}}(\lambda)= \frac{t\Gamma(2t)}{\lambda^t\;\sqrt{1-c^2}}\; \Ima\left(\er\left[\frac{e^{-\jcb(\lambda/\mathbb{S}_t)^t} }{\mathbb{S}_t^t}\right]\right)\nonumber\\
&=&\frac{t\Gamma(2t)}{\lambda^t\;\sin(\pi\al)}\; \er\left[\frac{e^{-\cos(\pi\al)(\lambda/\mathbb{S}_t)^t}}{\mathbb{S}_t^t} \; \sin\big(\sin(\pi\al)\; (\lambda /\mathbb{S}_t)^t\big)\right], \;\;\;\lambda \geq 0.
\end{eqnarray*}
3) This done as in proof of point 3) of Corollary \ref{crt1}.
\end{proofof}
\medskip

\medskip
\begin{proofof} {\bf Corollary \ref{corg}.} Equivalences $1) \Longleftrightarrow 2) \Longleftrightarrow 3)$ are due to point 2) of Proposition \ref{hcmlog} and to \eqref{party},  $2) \Longleftrightarrow 4)$ corresponds to  \eqref{factuel} $\Longleftrightarrow$ \eqref{fact}, and $1) \Longleftrightarrow 5)$ is due to point 2) of Proposition \ref{hcmlog} and to \eqref{lala}. For the last assertion, assume $\alpha \in (0,1/2)$, $0<t\leq  1-\alpha$,  $c=\cos(\pi \alpha)$ in \eqref{lala}, then consider the $\GGC$   r.v. $\mathbb{X}_{c,t}$, which is linked to $\tal^{1/t}$ by the expression
\begin{equation}\label{fact'}
f_{\tal^{1/t}}(x)=  \frac{t \sin(\pi \alpha)}{\pi \alpha}\,\frac{ x^{t-1}}{1+2 \cos(\pi \alpha) x^t + x^{2t}} = \frac{t \sin(\pi \alpha)}{\pi \alpha}\,  x^{t-1} \, \er[e^{-x \;\mathbb{X}_{c,t}}].
\end{equation}
Integrating both sides in  \eqref{fact'},  we see that
\begin{equation}\label{mom}
\er[(\mathbb{X}_{c,t})^{-t}]= \frac{\pi \alpha}{\sin(\pi \alpha) \Ga(t+1)}<\infty,
\end{equation}
which  by the procedure \eqref{biais},  allows to introduce the r.v. $\mathbb{X}_{c,t}^{\;[-t]}$, whose  distribution is also a $\GGC$, by \eqref{propel} and \eqref{gcgc}, we obtain the equivalences
$$\mathbb{X}_{c,t} \sim \GGC \Longleftrightarrow \er[e^{-x \;\mathbb{X}_{c,t}}] \in \hcm \Longleftrightarrow f_{\tal^{1/t}}(x) \in \hcm.$$
Using \eqref{ideng} and the effect of size-biasing in \eqref{biaiso}, we obtain
\begin{equation}\label{tog}
\mathbb{X}_{c,t}^{\;[-t]} \simdis \left(\frac{\ga_{2t}}{\dac}\right)^{[-t]}\simdis \frac{\ga_{t}}{\dac^{\;[t]}},
\end{equation}
hence $\mathbb{X}_{c,t}^{\;[-t]} \;$ is also a $\ga_t$-mixture.
Then, \eqref{fact} and \eqref{mom} yield
\begin{eqnarray*}
f_{\tal^{1/t}}(x) &=& \frac{t \sin(\pi \alpha)}{\pi \alpha}\,  x^{t-1} \, \er[e^{-x \;\mathbb{X}_{c,t}} ]
=\frac{t \sin(\pi \alpha)}{\pi \alpha}\,  x^{t-1} \, \er[(\mathbb{X}_{c,t})^t\, (\mathbb{X}_{c,t})^{-t} e^{-x \;\mathbb{X}_{c,t}} ]\\
&=&\frac{t \sin(\pi \alpha)}{\pi \alpha}\; \er[(\mathbb{X}_{c,t})^{-t}]\;  x^{t-1} \; \er[(\mathbb{X}_{c,t}^{\;[-t]})^t  e^{-x\; \mathbb{X}_{c,t}^{\;[-t]}} ]= \frac{x^{t-1}}{\Ga(t)} \, \er[(\mathbb{X}_{c,t}^{\;[-t]} )^t\,e^{-x\;\mathbb{X}_{c,t}^{\;[-t]} }]\\
&=& f_{\ga_t/\mathbb{X}_{c,t}^{\;[-t]} }(x),
\end{eqnarray*}
which, combined with \eqref{tog}, gives \eqref{idt}.
\end{proofof}

\section{Conclusion and perspectives}
In this paper, we showed that the Mittag-Leffler function (with eventually complex argument) are tightly linked to the stable distributions by various aspects: we  exhibited their non-trivial infinite divisibility, $\GGC$ and $\HCM$ properties, and we provided their explicit intervention in the distributional properties for the first passage time of the spectrally positive stable process. We also introduced new classes of $\HCM$ distributions, and gave a possible direction to solve the open question \eqref{conju} on the power of the positive stable r.v.'s.  Indeed, \eqref{idt} gives
$$\left(\frac{\sal}{\sal'}\right)^{\frac{\al}{t}} \simdis \frac{\ga_{t}}{(\mathbb{X}_{c,t})^{[-t]}}\simdis  \frac{\ga_t}{\ga_t'} \dac^{\;[t]}\sim \HCM, \;\;\mbox{if}\;\alpha\in (0,1/2], \;c=\cos(\pi\alpha) \;\mbox{and }\;  t=1-\alpha.$$
The latter indicates that two independent factorizations are feasible, namely,
$$\left(\sal\right)^{\frac{\al}{t}}\simdis \ga_{t}\; \mathbb{Z}_{c,t}, \quad \dac^{\;[t]} \simdis \frac{\mathbb{Z}_{c,t}}{\mathbb{Z}_{c,t}'}, \;\;\;\mathbb{Z}_{c,t},\;\mathbb{Z}_{c,t}'\;i.i.d.$$
In other terms the distribution of $\log \dac^{\;[t]}$ is symmetric. More investigation of the distribution of
$\dac^{\;[t]}$ is then necessary to solve the open question \eqref{conju}. 

\bigskip
\noindent{\bf Acknowledgment:}  We thank the two anonymous reviewers for their valuable comments and suggestions which substantially improved the presentation and the results of this article.

\end{document}